\documentclass[12pt,centertags]{amsart}
\usepackage{amsmath,amstext,amsthm,a4,amssymb,amscd}
\usepackage{newcent}
\usepackage[mathscr]{eucal}
\usepackage{hyperref}
\usepackage{mathrsfs}
\numberwithin{equation}{section}

\newcommand{\field}[1]{\mathbb{#1}}
\newcommand{\bZ}{\field{Z}}
\newcommand{\bR}{\field{R}}
\newcommand{\bC}{\field{C}}
\newcommand{\bN}{\field{N}}

\def\bJ{{\mathbf J}}

\newcommand{\cali}[1]{\mathscr{#1}}
\newcommand{\cC}{\cali{C}} 
\newcommand{\cO}{\cali{O}} 
\newcommand{\cH}{\cali{H}} 
 \newcommand{\cL}{\cali{L}}

\newcommand{\calig}[1]{\mathcal{#1}}

 \newcommand\mO{\calig{O}}
\newcommand\mQ{\calig{Q}} \newcommand\mR{\calig{R}}

\def\bE{{\mathbf E}}

\DeclareMathOperator{\End}{End} 
\DeclareMathOperator{\Ker}{Ker} 
 
  \DeclareMathOperator{\spec}{Spec}
\DeclareMathOperator{\Id}{Id} 
\DeclareMathOperator{\tr}{Tr} 
 
 \newcommand{\spin}{$\text{spin}^c$ }
\DeclareMathOperator{\db}{\overline\partial}

\newcommand{\norm}[1]{\lVert#1\rVert} \newcommand{\abs}[1]{\lvert#1\rvert}
\newcommand{\om}{\omega} 
\newcommand{\sbullet}{\scriptscriptstyle{\bullet}}

\newtheorem{thm}{Theorem}[section]

\theoremstyle{definition}

\theoremstyle{remark}
\newtheorem{rem}[thm]{Remark}
\newcommand{\be}{\begin{equation}}
\newcommand{\ee}{\end{equation}}
\newcommand{\wi}{\widetilde}
\newcommand{\var}{\varepsilon}
\newcommand{\ov}{\overline}
\newcommand{\comment}[1]{}

\begin{document}
\title[first coefficients of the asymptotic expansion of the Bergman kernel]
{The first coefficients of the asymptotic expansion of the Bergman kernel
of the \spin Dirac operator}
\date{\today}
 \author{Xiaonan Ma}
 \address{ Centre de Math\'ematiques Laurent Schwartz, 
 UMR 7640 du CNRS, \'Ecole Polytechnique, 91128 Palaiseau Cedex, France}
\email{ma@math.polytechnique.fr}


\author{George Marinescu}
\address{Fachbereich Mathematik, Johann Wolfgang Goethe-Universit\"at, 
Robert-Mayer-Stra{\ss}e 10, 60054, Frankfurt am Main, Germany}
\email{marinesc@math.uni-frankfurt.de}

\begin{abstract}
We establish the existence of the asymptotic expansion
of the Bergman kernel associated to the spin$^c$ Dirac operators
acting on high tensor powers of line bundles 
with non-degenerate  mixed curvature (negative and positive eigenvalues)
by extending \cite{DLM04a}. 
We compute the second coefficient $\pmb{b}_1$ 
in the asymptotic expansion using the method of \cite{MM04a}.
\end{abstract}

\maketitle

\setcounter{section}{-1}
\section{Introduction}

In \cite[Theorem 4.18$^\prime$]{DLM04a}, Dai, Liu and Ma established 
the full off-diagonal asymptotic expansion of the Bergman kernel  
of the \spin Dirac operator on high tensor powers of a line bundle
on a compact symplectic manifold.
This paper is a continuation of their work.

For some applications, it is necessary to 
compute the first two coefficients $\pmb{b}_0$, $\pmb{b}_1$ 
of the asymptotic expansion.
The approach of \cite{DLM04a} is to relate the heat kernel and the Bergman kernel
expansions. The computation of the coefficients of the Bergman kernel expansion is
done by using the corresponding coefficients of the heat kernel.
Thus, $\pmb{b}_0$ is calculated in \cite[Theorem 1.1]{DLM04a} and 
$\pmb{b}_1$ in \cite[Theorem 1.3]{DLM04a}, the latter only in the K\"ahler case.
These results brought a new proof of \cite{Catlin99, Lu00, Wang1, Ze98}.

Considering the symplectic case, the coefficient of $t^1$ 
in the Taylor expansion of the rescaled operator does not vanish,
thus it is complicate to compute $\pmb{b}_1$ in this way.
However, we developed in \cite[\S 1.5]{MM04a} a method of formal power series  
to compute the coefficients for the renormalized Bochner-Laplacian.
The main result of this paper is Theorem \ref{t0.2},
where we compute the coefficient $\pmb{b}_1$ 
in the asymptotic expansion
of the Bergman kernel associated to the spin$^c$ Dirac operators
by applying the method in \cite[\S 1.5, 2.3]{MM04a}.
Comparing with \cite[\S 2.3]{MM04a}, 
the contribution from the coefficient of $t^1$ 
in the Taylor expansion of the rescaled operator $L^t_2$ in \eqref{c30}
is quite complicate here. But after taking the trace in \eqref{c2},
we refind the Hermitian scalar curvature of \cite{D97}.
For more details on our approach we also refer the readers to 
our recent book \cite{MM05b}.

Another feature of the paper is that we do not suppose that
the almost complex structure polarizes the curvature of the line bundle,
that is we allow bundles with mixed curvature (negative and positive eigenvalues). 
In Section \ref{s5}, we explain that the arguments in 
\cite{DLM04a} still work well under this condition (cf. \eqref{diag1}).
Then we compute the coefficient $\pmb{b}_1$ in Section \ref{s2}.


\subsection*{Aknowledgements}
We thank Professor Johannes Sj\"ostrand for useful conversations
and the referee for  useful comments.

\section{Bergman kernel of the spin$^c$ Dirac operator}\label{s5}

This Section is organized as follows. In Section \ref{s5.1}, we
recall the Lichnerowicz formula for the \spin Dirac operator $D_p$.
As a consequence we exhibit the spectral gap 
for $D_p^2$ without assumption that the almost complex structure $J$
polarizes the symplectic form $\om$. This is done in Section \ref{s5.2}.
We explain in Section \ref{s5.3} 
the full off-diagonal asymptotic expansion for the Bergman kernel.
Then we show how to handle the operator 
$\ov{\partial}+\ov{\partial}^*$, which is the content of Section \ref{s5.4},
we explain also  its relation to the tangent Cauchy-Riemann complex 
in Section \ref{s5.5}.

\subsection{The \spin Dirac operator}\label{s5.1}

Let $(X,\om)$ be a compact connected symplectic manifold of real dimension
$2n$. Assume that there exists a Hermitian line bundle $L$ over
$X$ endowed with a Hermitian connection $\nabla^L$ with the
property that $\frac{\sqrt{-1}}{2\pi}R^L=\omega$, where
$R^L=(\nabla^L)^2$ is the curvature of $(L,\nabla^L)$. Let
$(E,h^E)$ be a Hermitian vector bundle  on $X$ with Hermitian
connection $\nabla^E$ and its curvature $R^E$.

Let $J$ be an almost complex structure on $TX$ and
$g^{TX}$ be a Riemannian metric on $X$ compatible with $J$,
i.e. $g^{TX}(\cdot, \cdot)= g^{TX}(J\cdot, J\cdot)$.
We designate by $\nabla ^{TX}$ the Levi-Civita connection on $(TX,
g^{TX})$ and by $R^{TX}$ and $r^X$ its curvature and scalar curvature, respectively.

The almost complex structure $J$ induces a splitting
$T X\otimes_\bR \bC=T^{(1,0)}X\oplus T^{(0,1)}X$,
where $T^{(1,0)}X$ and $T^{(0,1)}X$ are the eigenbundles of $J$ 
corresponding to the eigenvalues $\sqrt{-1}$ and $-\sqrt{-1}$ respectively.
Then $P^{(1,0)}= \frac{1}{2} (1-\sqrt{-1} J)$ 
is the projection from $TX\otimes_\bR  \bC$ onto $T^{(1,0)}X$.
For any $v\in TX\otimes_\bR \bC$ with decomposition
$v=v_{1,0}+v_{0,1} \in T^{(1,0)}X\oplus T^{(0,1)}X$, let
${\overline v^\ast_{1,0}}\in T^{*(0,1)}X$ be the metric dual of
$v_{1,0}$. Then $c(v)=\sqrt{2}({\overline v^\ast_{1,0}}\wedge-
i_{v_{0,1}})$ defines the Clifford action of $v$ on $\Lambda
(T^{*(0,1)}X)$, where $\wedge$ and $i$ denote the exterior and
interior product respectively.

Let $\nabla ^{T^{(1,0)}X} =P^{(1,0)} \nabla ^{TX}P^{(1,0)}$ 
be the connection on $T^{(1,0)}X$ induced by $\nabla ^{TX}$, 
with curvature $R^{T^{(1,0)}X}$. 
Let $\nabla^{\text{Cliff}}$ be the Clifford connection 
on $\Lambda (T^{*(0,1)}X)$ induced canonically by $\nabla^{TX}$
(cf. \cite[\S 2]{MM02})\footnote{In \cite[(2.3)]{MM02},
one missed a term '$+\frac{1}{2}\tr_{T^{(1,0)}X}\Gamma$' 
in the right hand side of the first line,
and the second line should be read as 
$'=d+\sum_{lm}\lbrace\big\langle\Gamma w_l,\overline{w}_m
\big\rangle\,\overline{w}^m\wedge\,i_{\overline{w}_l}+'$ } 
and let $R^{\text{Cliff}}$ be its curvature.
Let $\{e_i\}_{1\leqslant i\leqslant 2n}$ be an orthonormal basis of $TX$.
From the definitions, we get (cf. also \cite[(0.12)]{MM04a})
\begin{equation}\label{c20}
\begin{split}
&R^{T^{(1,0)}X} = P^{1,0} \Big[R^{TX}
-\frac{1}{4} (\nabla ^XJ)\wedge (\nabla ^XJ)\Big]P^{1,0}\,,\\
&  R^{\text{Cliff}}
=\frac{1}{4} \left\langle R^{TX} e_l,e_m  \right\rangle 
\, c(e_l)\, c(e_m) + \frac{1}{2} \tr\left[R^{T^{(1,0)}X}\right].
\end{split}
\end{equation}
Let $\nabla^{E_p}$ be the connection on
$E_p:=\Lambda (T^{*(0,1)}X)\otimes L^p\otimes E$
induced by $\nabla^{\text{Cliff}}$, $\nabla^L$ and  $\nabla^E$.
The \spin Dirac operator $D_p$ is defined using  
$\nabla^{E_p}$ and acts on
$\Omega^{0,{\scriptscriptstyle{\bullet}}}(X,L^p\otimes E)
=\bigoplus_{q=0}^n\Omega^{0,q}(X,L^p\otimes E)$, the direct sum of
spaces of $(0,q)$--forms with values in $L^p\otimes E$. We refer the reader to 
\cite[(2.6)]{MM02} for the precise definition. We content ourselves to describe
the operator $D^2_p$ in the sequel.

Let $\langle\,\cdot\,,\,\cdot\,\rangle_{E_p}$ be the metric on $E_p$ induced by
 $g^{TX}$, $h^L$ and $h^E$.
Let $dv_X$ be the Riemannian volume form of $(TX, g^{TX})$.
The $L^2$--scalar product on $\Omega^{0,\scriptscriptstyle{\bullet}}(X,L^p\otimes E)$, the space of smooth sections of
$E_p$, is given by
\begin{equation}\label{b3}
\langle s_1,s_2 \rangle =\int_X\langle s_1(x),
s_2(x)\rangle_{E_p}\,dv_X(x)\,.
\end{equation}
We denote the corresponding norm with $\norm{\cdot}_{L^ 2}$.

Let $\left(\nabla^{E_p}\right)^{\ast}$ be the formal adjoint of 
$\nabla^{E_p}$ with respect to \eqref{b3}.
Set 
\begin{equation}\label{diag0}
\mathbf{c}(R)= \frac{1}{2}  
\left(R^E+\tfrac{1}{2}\tr[R^{T^{(1,0)}X}]\right)(e_l,e_m)
\,c(e_l)\,c(e_m).
\end{equation}
Then the Lichnerowicz formula \cite[Theorem 3.52]{BeGeVe}
(cf. \cite[Theorem 2.2]{MM02}) for $D_p^2$ is
\begin{equation}\label{Lich}
D^2_p=\left(\nabla^{E_p}\right)^{\ast}\,\nabla^{E_p}
+\frac{1}{4}r^X 
+  \frac{1}{2} p\, R^L(e_l,e_m) \,c(e_l)\,c(e_m) + \mathbf{c}(R).
\end{equation}

The Bergman kernel  $P_p(x,x')$,  $(x,x'\in X)$, 
is the smooth kernel with respect to $dv_X(x')$ of the orthogonal projection $P_p$ from
$\Omega^{0,\sbullet}(X,L^p\otimes E)$
on $\Ker D_p$\,. 
Then $P_p(x,x)$ is an element of $\End (\Lambda (T^{*(0,1)}X)\otimes E)_x$\,.

\subsection{Spectral gap of the spin$^c$ Dirac operator}\label{s5.2}
We choose the almost complex structure $J$ such that $\om$ is $J$-invariant, 
i.e. $\om(\cdot,\cdot)=\om(J\cdot,J\cdot)$. 
But we do not suppose that  $\om(\cdot,J\cdot)$ is positive 
in Sections \ref{s5.2}-\ref{s5.4}. This is the difference comparing with 
the assumption in \cite{MM02, DLM04a, MM04a}.

Let $\bJ:TX\longrightarrow TX$ be the skew--adjoint linear
map which satisfies the relation
\be \label{a0.1}
\om(u,v)=g^{TX}(\bJ u,v)
\ee
 for $u,v \in TX$. Then $J$ commutes with $\bJ$. 
Thus $\bJ \in \End(T^{(1,0)}X)$, and for any $x\in X$, 
we can diagonalize $\bJ_{x}$, i.e. find an orthonormal
 basis $\{w_j\}_{j=1}^n$ of $T^{(1,0)}X$ such that 
$\bJ_{x} w_j =\frac{\sqrt{-1}}{2\pi}a_j(x) w_j$ with $a_j(x)\in \bR$.
As $\om$ is non-degenerate, the number of negative
eigenvalue of $\bJ_{x}\in \End(T^{(1,0)}_x X)$ does not depend on $x$,
and we denote it by $q$ (In \cite{MM02, DLM04a, MM04a}, we suppose that $q=0$).
From now on, we assume that 
\begin{align}\label{diag1}
 \bJ_{x} w_j =\frac{\sqrt{-1}}{2\pi}a_j\,  w_j, \quad
a_j(x)<0 \, \,\mathrm{for}\,\, j\leqslant q \,\,\mathrm{and}
\, \,a_j(x)>0\,\, \mathrm{for}\,\, j> q.
\end{align}
Then the vectors $\{w_j\}_{j=1}^q$ span a sub-bundle $W$ of $T^{(1,0)} X$.
Set
\begin{equation}\label{diag2}
\begin{split}
&\om_d(x) =- \sum_{j=1}^n a_j \ov{w}^j \wedge i_{\ov{w}_j}
 + \sum_{j=1}^q a_j
=\sum_{j=1}^q a_j i_{\ov{w}_j}\wedge \ov{w}^j    
- \sum_{j=q+1}^n a_j \ov{w}^j \wedge i_{\ov{w}_j} \,,\\
&\tau(x)= \pi \tr|_{TX}  (-\bJ^2)^{1/2} = \sum_{j=1}^n |a_j|
=-\sum_{j=1}^q a_j+ \sum_{j=q+1}^n a_j\,,\\
&\mu_0= \inf_{ x\in X, \, j}|a_j(x)| \,.
\end{split}
\end{equation}
Then we have 
\begin{align}\label{diag3}
  \frac{1}{2} R^L(e_l,e_m)\,c(e_l)\,c(e_m)
= -2 \om_d -\tau.
\end{align}

The following result extends \cite[Theorem 2.5]{MM02} 
to the current situation. We denote by 
$\Omega^{\neq q}(X,L^p\otimes E)
=\bigoplus_{k\neq q}\Omega^{0,k}(X,L^p\otimes E)$.
\begin{thm}\label{main}
There exists $C>0$ such that for any $p\in \bN$  
\begin{equation}\label{main1}
\norm{D_{p}s}^2_{L^2}\geqslant(2p\mu_0-C)\norm{s}^2_{L^2}\,,\quad\text{for $s\in\Omega^{\neq q}(X,L^p\otimes E)$}.
\end{equation}
\end{thm}
\begin{proof}
By \eqref{Lich} and \eqref{diag3}, 
for $s\in\Omega^{0,\scriptscriptstyle{\bullet}}(X,L^p\otimes E)$\,,
\begin{equation}\label{lich1}
\norm{D_{p}s}^2_{L^2}=\lbrace\norm{\nabla^{\Lambda^{0,\scriptscriptstyle{\bullet}}\otimes L^p\otimes E}s}^2_{L^2}-p\langle\tau(x)s,s\rangle\rbrace
-2p \left\langle\om_{d}s,s\right\rangle+
\left\langle\left(\tfrac{1}{4}r^X +\mathbf{c}(R)\right)s,s\right\rangle\,.
\end{equation}
We consider now $s \in \cC^\infty(X,L^p\otimes E')$,
 where $E^\prime=E\otimes\Lambda (T^{*(0,1)}X)$.  
By \cite[Corollary 2.4]{MM02} which is a direct consequence 
of the Lichnerowicz formula 
(cf. also \cite[Theorem 1]{GuU88}, \cite[Theorem 2.1]{BorU96}, 
\cite[Theorem 4.4]{Bra99}),  there exists $C>0$ such that for any $p>0$, 
$s \in \cC^\infty(X,L^p\otimes E')$, we have 
\begin{equation}\label{est1prim}
\big\|\nabla^{L^p\otimes E^\prime}s\big\|^2_{L^2}
-p\left\langle\tau(x)s,s\right\rangle
\geqslant -C\norm{s}^2_{L^2}\,.
\end{equation}
If $s\in\Omega^{\neq q}(X,L^p\otimes E)$, 
the second term of (\ref{lich1}), 
$-2p\left\langle\om_{d}s,s\right\rangle$ is bounded below by
$2p \mu_0 \norm{s}^2_{L^2}$, while the third term of (\ref{lich1}) is $\cO(\norm{s}^2_{L^2})$.
The proof of \eqref{main1} is completed.
\end{proof}

Set  
 \begin{equation}
 \begin{split}\label{2.64}
D_p^{+}&=D_p|_{\Omega^{0,{\text{even}}}}, \quad
 D_p^{-}=D_p|_{\Omega^{0,{\text{odd}}}}, \\
o_q&= - \quad \mathrm{if} \,\, q \,\, \mbox{is even};\quad 
  o_q  =+\quad \mathrm{if} \,\, q\,\,   \mbox{is odd}\,. 
\end{split}
\end{equation}

By using the trick of the proof of Mckean-Singer formula
and Theorem \ref{main}, the same proof as in \cite[\S 3]{MM02}
gives the following extension of \cite[Theorem 1.1]{MM02}.
For any operator $A$, we denote by $\spec (A)$ the spectrum of $A$. 
\begin{thm}\label{specDirac}
There exists $C>0$ such that for $p\in \bN$, 
 \begin{align}\label{diag5}
\spec (D^2_p) \subset \{0\}\cup [2p\mu_0-C,+\infty[. 
\end{align}
For $p$ large enough,
we have 
\begin{equation}\label{vanish}
\Ker D_p^{o_q}=\{0\}.
\end{equation}
\end{thm}

\subsection{Off-diagonal asymptotic expansion of Bergman kernel}\label{s5.3}

The existence of the spectral gap expressed in Theorem \ref{specDirac}
allows us to obtain immediately as in \cite[Prop. 4.1]{DLM04a} the far
off-diagonal behavior of the Bergman kernel.
Namely, for any $l,m\in \bN$ and $\var>0$, there exists $C_{l,m,\var}>0$ 
such that for $p\geqslant 1$, $x,x'\in X$, $d(x,x')> \var$, 
\begin{align}\label{1c3}
&|P_p(x,x')|_{\cC^m(X\times X)} \leqslant C_{l,m,\var}\, p^{-l}.
\end{align}

We denote by
$ I_{\det(\ov{W}^*)\otimes E}$  the orthogonal
projection from $\bE:=\Lambda (T^{*(0,1)}X)\otimes E$
 onto $\det(\ov{W}^*)\otimes E$.
Let $\pi : TX\times_{X} TX \to X$ be the natural projection from the
fiberwise product of $TX$ on $X$.
  Let $\nabla ^{\End (\bE)}$ be the connection on
$\End (\Lambda (T^{*(0,1)}X)\otimes E)$ induced by
$\nabla ^{\text{Cliff}}$ and $\nabla ^E$.

For $x_0\in X$, we identify $L_Z$, $E_Z$ and $(E_p)_Z$
for $Z\in B^{T_{x_0}X}(0,\var)$ to $L_{x_0}, E_{x_0}$ and $(E_p)_{x_0}$
by parallel transport with respect to the connections
$\nabla ^L, \nabla ^E$ and
 $\nabla^{E_p}$ along the curve $\gamma_Z :[0,1]\ni u \to \exp^X_{x_0} (uZ)$.
Under this identification and \eqref{1c3}, we will view 
$P_p(x,x')$ as a smooth section 
$P_{p,x_0}(Z,Z')$, $(Z,Z'\in B^{T_{x_0}X}(0,\var))$,
 of $\pi ^* (\End (\Lambda (T^{*(0,1)}X)\otimes E))$ on  $TX\times_{X} TX$.
 And $\nabla ^{\End (\bE)}$
induces naturally a $\cC^m$-norm for the parameter $x_0\in X$.

 Let $dv_{TX}$ be the Riemannian volume form on
$(T_{x_0}X, g^{T_{x_0}X})$.
Let $\kappa (Z)$ be the smooth positive function defined by the equation
\be\label{c22}
dv_{X}(Z) = \kappa (Z) dv_{TX}(Z),
\ee
with $\kappa(0)=1$.
We denote by $\det_{\bC}$ for the determinant function on the complex
 bundle $T^{(1,0)}X$, and $|\bJ_{x_0}|=(-\bJ^2_{x_0})^{1/2}$. 
 Denote by  $\nabla_U$ the ordinary differentiation
 operator on $T_{x_0}X$ in the direction $U$.
On $T_{x_0}X\simeq \bR^{2n}$, set
\begin{align} \label{3ue63}
L^0_{2,\bC}=- \sum_j 
\Big(\nabla _{e_j} +\frac{1}{2}R^L_{x_0}(Z, e_j)\Big)^2   -\tau_{x_0}.
\end{align}
 Let $P(Z,Z')$ be the Bergman kernel of $L^0_{2,\bC}$,
 i.e. the smooth kernel of the orthogonal projection from 
$L^2 (\bR^{2n},\bC)$ onto $\Ker L^0_{2,\bC}$. 
Then for $Z,Z'\in T_{x_0}X$,
\begin{align}\label{3ue62}
P(Z,Z')  ={\det}_{\bC} (|\bJ_{x_0}|)
\exp\Big (- \frac{\pi}{2} \left \langle |\bJ_{x_0}|(Z-Z'),(Z-Z') \right \rangle
-\pi \sqrt{-1} \left \langle \bJ_{x_0} Z,Z' \right \rangle\Big ).
\end{align}

The main result of this part is the following extension of 
\cite[Theorem 4.18$^\prime$]{DLM04a}.
\begin{thm} \label{tue17}
There exist polynomials $J_{r}(Z,Z^{\prime})\in
\End (\Lambda (T^{*(0,1)}X) \otimes E)_{x_0}$ $(x_0\in X)$,
in $Z,Z^{\prime}$ with the same parity as $r$
and with $\deg J_{r}\leqslant 3r$, whose coefficients
are polynomials in $R^{TX}$, $R^{T^{(1,0)}X}$,
$R^E$ {\rm (}and $R^L${\rm )}  and their derivatives of order
$\leqslant 2r-1$  {\rm (}resp. $\leqslant2r${\rm )} and reciprocals
of linear combinations of eigenvalues of $\bJ$  at $x_0$, 
 such that by setting
\begin{equation}\label{a0.7}
P^{(r)}_{x_0}(Z,Z^{\prime})=J_{r}(Z,Z^{\prime})P(Z,Z^{\prime}),
\quad J_{0}(Z,Z^{\prime})
= {\det}_\bC (|\bJ_{x_0}|)I_{\det(\ov{W}^*)\otimes E}\,\, ,
\end{equation}
the following statement holds:
There exist $C''>0$ such that 
for any $k,m,m'\in \bN$, there exist $N\in \bN$ and $C>0$ with 
\begin{multline}\label{aue66}
\left |\frac{\partial^{|\alpha|+|\alpha'|}}
{\partial Z^{\alpha} {\partial Z'}^{\alpha'}}
\left (\frac{1}{p^n}  P_p(Z,Z')
-\sum_{r=0}^k  P^{(r)} (\sqrt{p} Z,\sqrt{p} Z')\kappa ^{-1}(Z')
p^{-r/2}\right )\right |_{\cC^{m'}(X)}\\
\leqslant C  p^{-(k+1-m)/2}  (1+|\sqrt{p} Z|+|\sqrt{p} Z'|)^N
\exp (- \sqrt{C''\mu_0 } \sqrt{p} |Z-Z'|)+ \cO(p^{-\infty}).
\end{multline}
for any $\alpha, \alpha' \in \bN^{n}$, with $|\alpha|+|\alpha'|\leqslant m$, 
any $Z,Z'\in T_{x_0}X$ with $|Z|, |Z'|\leqslant  \var$ and any $x_0\in X$, $p\geqslant 1$.
\end{thm}
Here $\cC^{m'}(X)$ is the  $\cC^{m'}$-norm for the parameter $x_0\in X$.
We say that a term $T=\cO(p^{-\infty})$ if for any $l,l_1\in \bN$, 
there exists $C_{l,l_1}>0$ such that the $\cC^{l_1}$-norm of $T$ is dominated 
by $C_{l,l_1} p^{-l}$.
\begin{proof} By using Theorems \ref{main}, \ref{specDirac}, 
we know the arguments in \cite[\S 4.1-4.3]{DLM04a} go through without
any change.  By \eqref{Lich}, \eqref{diag3},  
as in \cite[(4.105)]{DLM04a}, the corresponding limit operator here is still
\begin{align}\label{ue54}
&L^0_2 = L^0_{2,\bC}-2 \omega_{d,x_0}.
\end{align}
 Let $e ^{-uL^0_{2,\bC}} (Z,Z')$, $e ^{-uL^0_2} (Z,Z')$ 
be the smooth kernels of $e ^{-uL^0_{2,\bC}}$, $e ^{-uL^0_2}$
with respect to $dv_{TX}(Z')$. Now from (\ref{ue54})
(cf. \cite[(6.37), (6.38)]{B90d}), 
we need to replace \cite[(4.106)]{DLM04a} by the following equations
\begin{equation}\label{ue56}
\begin{split}
&e ^{-uL^0_{2,\bC}} (Z,Z')
= {\det}_{\bC} \Big (\frac{|\bJ_{x_0}|}{1-e^{-4\pi u|\bJ_{x_0}|}}\Big )
\exp  \Big (-\frac{1}{2} \left \langle
\frac{\pi |\bJ_{x_0}|}{\tanh (2\pi u|\bJ_{x_0}|)} Z,Z \right \rangle\\
&\hspace*{6mm} -\frac{1}{2} \left \langle \frac{\pi |\bJ_{x_0}|}
{\tanh (2\pi u|\bJ_{x_0}|)} Z',Z' \right \rangle
+ \left \langle \frac{\pi |\bJ_{x_0}|}{\sinh (2\pi u|\bJ_{x_0}|)}
e ^{- 2\pi \sqrt{-1}u \bJ_{x_0}} Z,Z' \right \rangle\Big ),\\
&e ^{-uL^0_{2}} (Z,Z') = e ^{-uL^0_{2,\bC}} (Z,Z') e ^{2u \omega_{d,x_0}}.
\end{split}
\end{equation}
Observe that for $\omega_{d,x_0}\in \End (\Lambda (T^{*(0,1)}X))_{x_0}$, 
by \eqref{diag2},
\begin{equation} \label{3ue56} 
\begin{split}
&\Ker \omega_{d,x_0} = \det (\ov{W}^*)_{x_0},\\
&\omega_{d,x_0} \leqslant -\mu_0 \quad   \mbox{on} \, \,
(\det (\ov{W}^*)_{x_0})^\bot
\mbox{= the orthogonal complement of}\,  \det (\ov{W}^*). 
\end{split}
\end{equation}
Thus from \cite[\S 4.4, 4.5]{DLM04a}, we get Theorem \ref{tue17}.
\end{proof}
If we take $Z,Z'=0$ in Theorem \ref{tue17} we infer the following result
by the same argument in \cite[\S 4.5]{DLM04a}.
\begin{thm}\label{at0.1} There exist smooth coefficients
$\pmb{b}_r(x)\in \End (\Lambda (T^{*(0,1)}X)\otimes E)_x$ which
are polynomials in $R^{TX}$, $R^{T^{(1,0)}X}$,
$R^E$ {\rm (}and $R^L${\rm )} and their derivatives of order
$\leqslant 2r-1$  {\rm (}resp. $\leqslant2r${\rm )} and reciprocals
of linear combinations of eigenvalues of ${\bf J}$ at $x$, such that
 \begin{equation}\label{b0}
 \pmb{b}_0={\det}_\bC (|\bJ|) I_{\det (\ov{W}^*)\otimes E}
 \end{equation}
and for any $k,l\in \bN$, there exists
$C_{k,l}>0$ with
\begin{align}\label{a0.3}
&\Big |P_p(x,x)
- \sum_{r=0}^{k} \pmb{b}_r(x) p^{n-r} \Big |_{\cC^l} \leqslant C_{k,l} p^{n-k-1}.
\end{align}
for any $x\in X$ and $p\in \bN$\,.
Moreover,
the expansion is uniform
in that for any  $k,l\in \bN$, there is an integer $s$ such that
if all data {\rm(}$g^{TX}$, $h^L$, $\nabla ^L$, $h^E$, $\nabla ^E${\rm)}
run over a set which are bounded in $\cC^s$ and with $g^{TX}$
 bounded below, there exists the constant  $C_{k,\,l}$
independent of $g^{TX}$,
and the $\cC^l$-norm in \eqref{a0.3} includes also the
derivatives on the parameters\,.
\end{thm}
\subsection{Holomorphic case revisited}\label{s5.4}

In this Section, we suppose that $(X, J)$ 
is a complex manifold with complex structure $J$ 
and $E,L$ are holomorphic vector bundles on $X$. 
We assume that $\nabla ^E$, $\nabla ^L$ are the holomorphic Hermitian 
(i.e. Chern) connections on $(E,h^E)$, $(L,h^L)$ and moreover, 
$\om:=\frac{\sqrt{-1}}{2 \pi} R^L$ defines a symplectic form on $X$. 
Therefore the signature of the curvature $\frac{\sqrt{-1}}{2 \pi} R^L$
(i.e. number of negative and positive eigenvalues) with respect to any Riemannian
metric compatible with $J$ will be the same.
Let $g^{TX}$ be any Riemannian metric on $TX$ compatible with $J$. Since $g^{TX}$ is 
not necessarily K\"ahler, $\bJ\neq J$ in \eqref{a0.1} in general.
Recall the number $q$ is defined by \eqref{diag1}, i.e. is the number of negative
eigenvalues of $\om$. Set 
 \begin{align} \label{f11}
\Theta(X,Y)=  g^{TX}(JX,Y).
\end{align}
Then the 2-form $\Theta$ need not be closed 
(the convention here is different to 
\cite[(2.1)]{B87c} by a factor $-1$). 

Let $\overline{\partial} ^{L^p\otimes E,*}$ be the formal adjoint of
the Dolbeault operator $\overline{\partial} ^{L^p\otimes E}$ 
on the Dolbeault complex
 $\Omega ^{0,\sbullet}(X, L^p\otimes E)$ with the scalar product 
induced by $g^{TX}$, $h^L$, $h^E$ as in \eqref{b3}.
Set 
\begin{equation}\label{Dp}
D_p = \sqrt{2}\big(\, \overline{\partial} ^{L^p\otimes E}
+ \,\overline{\partial} ^{L^p\otimes E,*}\big)\,.
\end{equation}
We denote by $\Box^{L^p\otimes E}= \overline{\partial} ^{L^p\otimes E}\,\overline{\partial} ^{L^p\otimes E,*}
+\,\overline{\partial} ^{L^p\otimes E,*}\,\overline{\partial} ^{L^p\otimes E}$
the Kodaira-Laplacian. Then $D^2_p=2\Box^{L^p\otimes E}$ 
it is twice the Kodaira-Laplacian and preserves the $\bZ$-grading of $\Omega ^{0,\sbullet}(X, L^p\otimes E)$.
By Hodge theory, we know that for any $k, p\in \bN$, 
\begin{equation} \label{f12}
\Ker D_p|_{\Omega ^{0,k}}  =\Ker D_p^2|_{\Omega ^{0,k}}  
\simeq H^{0,k} (X,L^p\otimes E),
\end{equation}
where $H^{0,\sbullet} (X,L^p\otimes E)$ is the Dolbeault cohomology.
Here $D_{p}$ is not a spin$^c$ Dirac operator on 
$\Omega ^{0,\sbullet}(X, L^p\otimes E)$, and $D^2_p$ is not 
a renormalized Bochner--Laplacian as in \cite{MM04a}, so we cannot apply directly
Theorems \ref{main} and \ref{specDirac}. 
Now we explain how to recover the conclusions of these theorems in the case of $D_p$\,.
The first step is to exhibit the spectral gap.
\begin{thm}\label{nonkahler1}
 The statements of Theorems \ref{main} and \ref{specDirac} still hold for the operator $D_{p}$ defined by \eqref{Dp}.
In particular, for $p$ large enough,
\begin{equation} \label{af12}
H^{0,k} (X,L^p\otimes E)=0 \quad \mbox{for} \, \, k\neq q.
\end{equation}
\end{thm}
\begin{proof}
As we will use \cite[Theorem 2.3]{B89a}
to study the Bergman kernel in the sequel, we prove Theorem \ref{nonkahler1}
by explaining \cite[Theorem 2.3]{B89a}.

Let $S^{-B}$ denote the 1-form with values in the antisymmetric
 elements of $\End(TX)$ which satisfies
\begin{equation}\label{f13}
\langle S^{-B}(U)V,W  \rangle = - \frac{\sqrt{-1}}{2} 
\Big( (\partial- \overline{\partial} )\Theta\Big)(U,V,W)\,,\quad\text{for $U,V,W \in TX$}.
\end{equation}
The Bismut connection $\nabla ^{-B}$ on $TX$ is defined by 
\begin{align}\label{f14}
\nabla ^{-B} = \nabla ^{TX} +  S^{-B}.
\end{align}
Then by \cite[Prop. 2.5]{B89a}, $\nabla ^{-B}$ preserves the metric $g^{TX}$
 and the complex structure of $TX$. 
Let $\nabla ^{\det}$ be the holomorphic Hermitian connection on 
$\det (T^{(1,0)}X)$ whose curvature is denoted $R^{\det}$.
Then these two connections induce naturally 
an unique connection on $\Lambda (T^{*(0,1)}X)$
which preserves its $\bZ$-grading, and with the connections 
$\nabla ^L, \nabla ^E$, we get a connection $\nabla ^{-B,E_p}$ on 
$\Lambda (T^{*(0,1)}X)\otimes L^p\otimes E$. 
Let $(\nabla^{-B,E_p})^*$ be the formal adjoint of $\nabla^{-B,E_p}$.
Let $C(TX)$ be the Clifford bundle of $TX$.
We define a map $^c: \Lambda (T^*X)  \to C(TX)$, by sending 
$e ^{i_1}\wedge\cdots \wedge e ^{i_j}$ to $c(e_{i_1})\cdots c(e_{i_j})$
for $i_1< \cdots < i_j$. 
For $B\in \Lambda^3(T^*X)$, set
$|B|^2= \sum_{i<j<k}|B(e_i,e_j,e_k)|^2$.
Then we can formulate  \cite[Theorem 2.3]{B89a} as following:
\begin{multline}\label{f16}
D_p^2= (\nabla^{-B,E_p})^*  \nabla^{-B,E_p} + \frac{r^X}{4} + {^c(R^E +pR^L 
+ \frac{1}{2} R^{\det})} \\
+ \frac{\sqrt{-1}}{2} {^c(\overline{\partial}\partial \Theta) }
- \frac{1}{8} |( \partial- \overline{\partial} )\Theta|^2.
\end{multline}
(\eqref{f16} can be seen as a Bochner-Kodaira-Nakano type formula.)
By using \eqref{diag3}, \eqref{est1prim} and \eqref{f16}, 
as in Theorems \ref{main} and \ref{specDirac}, we see that the conclusions
of these theorems still hold for the operator $D_p$ defined in \eqref{Dp}.
In particular \eqref{main1} holds.
Now from \eqref{main1} and  \eqref{f12}, we get \eqref{af12}.
\end{proof}

\begin{rem}\label{kahler2} The vanishing result \eqref{af12} 
is  Andreotti-Grauert's coarse vanishing theorem \cite[\S 23]{AdGr62}
(where it is proved by using the cohomology finiteness theorem for the disc
bundle of $L^*$). 
It can also be deduced, as shown by Griffiths \cite[p.\,432]{Gri66},
from the usual Bochner-Kodaira-Nakano formula
\cite{dem2,Gri66,Oh82}. The latter implies the following inequality for
$u\in\Omega^{m,k}(X,L^p\otimes E)$:
\begin{multline}\label{bkn}
\frac{3}{2}\,\Big(\norm{\db^{L^p\otimes E}u}^2+\norm{\db^{L^p\otimes E,*}u}^2\Big)
\geqslant\int_X\big\langle[\sqrt{-1}\,(pR^L+R^E), \Lambda]u,u\big\rangle\,dv_X
\\-\frac{1}{2}\,\big(\norm{Tu}^2+\norm{T^*u}^2+\norm{\overline Tu}^2+
\norm{\overline T^*u}^2\big)
\end{multline} 
where $\Lambda=i(\Theta)$ denotes the interior product with $\Theta$ and 
$T=[\Lambda,\partial\Theta]$ is the torsion of the metric $g^{TX}$.
We have pointwise
\begin{equation}\label{bkn2}
\big\langle[\sqrt{-1}\,R^L, \Lambda]u,u\big\rangle\geqslant
(\lambda_1+\dotsc+\lambda_k-\lambda_{n-m+1}-\dotsc-\lambda_n)\abs{u}^2,
\end{equation}
where $\lambda_1\leqslant\lambda_2\leqslant\dotsc
\leqslant\lambda_n$ are the eigenvalues of $\sqrt{-1}\,R^L$ with respect to $\Theta$.
If $m=0$ the right-hand side becomes $(-\lambda_{k+1}-\dotsc-\lambda_n)\abs{u}^2$.
As in \cite[Lemma\,4.3]{Oh82} we can {\em restrict} ourselves to those metrics $g^{TX}$ such that
the negative eigenvalues $\lambda_1,\dotsc,\lambda_q$ are very large and the
positive ones $\lambda_{q+1},\dotsc,\lambda_n$ are very small in absolute value.
Therefore, for $k<q$ there exists a constant $\mu_1>0$ such that
$-\lambda_{k+1}-\dotsc-\lambda_n\geqslant \mu_1$ on $X$. By \eqref{bkn}
we obtain \eqref{main1} for $u\in\Omega^{0,<q}(X,L^p\otimes E)=\oplus_{k<q}\Omega^{0,k}(X,L^p\otimes E)$.

In order to consider the case $k>q$, we apply again \eqref{bkn} and \eqref{bkn2} for 
$(n,q)$ $L^p\otimes K^{*}_X\otimes E$--valued forms, which involves
another change of metric, for which $\lambda_1,\dotsc,\lambda_q$ are small and
$\lambda_{q+1},\dotsc,\lambda_n$ are large in absolute value.
Thus we get \eqref{main1} also for $u\in\Omega^{0,>q}(X,L^p\otimes E)=
\oplus_{k>q}\Omega^{0,k}(X,L^p\otimes E)$,
but for yet another class of metrics $g^{TX}$.
Of course, the estimates just obtained entail immediately \eqref{af12}. 

We see however that by using \eqref{bkn} the essential estimate
\eqref{main1} for a fixed metric $g^{TX}$ seems out of reach, as well as 
the existence of the spectral gap \eqref{diag5}.
\end{rem}

By Theorem \ref{nonkahler1} the kernel of $D^2_p$ is concentrated in degree $q$. We consider
thus the Bergman kernel of $D^2_p$ in this particular degree.
Let $P^{0,q}_p(x,x')$ be the smooth kernel with respect to $dv_X(x')$ of the orthogonal 
projection from 
$\Omega^{0,q}(X,L^p\otimes E)$ on $\Ker D_p^2$\,.

\begin{thm}\label{nonkahler}
The Bergman kernel $P^{0,q}_p(x,x')$
has a full off--diagonal asymptotic expansion analogous to \eqref{aue66}
with $J_0=\det_\bC(|\bJ|) I_{\det(\ov{W}^*)\otimes E}$ 
as $p\to\infty$\,.
\end{thm}
\begin{proof}
We use now the connection $\nabla ^{-B,E_p}$ instead of $\nabla ^{E_p}$ 
in \cite[\S 3]{DLM04a}. Then by \eqref{main1} and \eqref{f16}, everything
goes through perfectly well and as in \cite[Theorem 4.18]{DLM04a},
so we can directly apply the result from \cite{DLM04a} to get the 
{\em full off-diagonal} asymptotic expansion of the Bergman kernel. 
As the above construction preserves the $\bZ$-grading on 
$\Omega ^{0,\sbullet}(X, L^p\otimes E)$, we can directly work
on $\Omega^{0,q}(X,L^p\otimes E)$.
\end{proof}
\begin{rem}\label{kahler4} From the arguments here and 
\cite{MM02}, \cite[\S 3.5]{MM04a}, 
we get naturally the covering version of Sections \ref{s5.2}-\ref{s5.4}.  
\end{rem}

\subsection{Relation to the tangential Cauchy-Riemann complex}\label{s5.5}

If $q=0$, i.e. $\bJ$ has only positive eigenvalues,
Theorem \ref{nonkahler} boils down to \cite[Theorem 3.9]{MM04a}.
Theorem \ref{nonkahler} for $x=x'$ is first due to Zelditch \cite{Ze98} and
Catlin \cite{Catlin99} and is based on the Boutet de Monvel-Sj\"ostrand 
parametrix \cite{BouSj76} for the Szeg\"o projector on CR functions 
on the boundary of the ``Grauert tube'' associated to $L$.
For general $q\neq0$, Berman and Sj\"ostrand \cite{BerSj05}
recently studied the asymptotic expansion $P_p(x,x')$, too.
They use an approach of Melin-Sj\"ostrand originating in the theory of
Fourier integral operators with complex phase.

 
In this section we briefly discuss the link between our analysis for $q\neq0$
and the the kernel of the Szeg\"o projector on $(0,k)$ forms on the boundary of the Grauert tube.
We use the notations and assumptions from Section \ref{s5.4}.

Let $Y= \{u\in L^*, |u|_{h^{L^*}}=1\}$ be the unit circle bundle in $L^*$.
$Y$ is a real hypersurface in the complex manifold $L^*$ which the boundary of the disc bundle $D= \{u\in L^*, |u|_{h^{L^*}}<1\}$,
with defining function $\varrho=|u|_{h^{L^*}}-1$. The Levi form of $\varrho$
restricted to the complex tangent plane of $Y$ coincides with the pull-back
of $\om$ through the canonical projection $\pi:Y\to X$. Hence it has $q$ negative
and $n-q$ positive eigenvalues.
We denote by $T^{*(0,1)}(Y)=T^{*(0,1)}L^*\cap(T^*Y\otimes_\bR \bC)$ 
the bundle of $(0,1)$-forms tangential to $Y$ and by 
$\Omega^{0,k}(Y)$ the space of smooth sections of $\Lambda^{k}(T^{*(0,1)}Y)$.
The $\db$ operator on the ambient manifold $L^*$ induces as usual a tangential
Cauchy-Riemann complex on the hypersurface $Y$ \cite{KoRo65, AnHi1, AnHi2}.
\begin{equation}\label{CRtangential}
0\longrightarrow\Omega^{0,0}(Y)\stackrel{\db_b}{\longrightarrow}\Omega^{0,1}(Y)
\stackrel{\db_b}{\longrightarrow}\dotsm\stackrel{\db_b}{\longrightarrow}\Omega^{0,n}(Y)
\longrightarrow0.
\end{equation}
The $\db_b$ operator commutes with the action of $S^1$ on $Y$.

The connection $\nabla ^L$ on  $L$ induces a connection 
on the $S^1$-principal bundle $\pi : Y\to X$, 
and let $T^H Y \subset TY$ be the corresponding horizontal bundle.
Let us introduce the Riemannian metric $g^{TY}=\pi^*(g^{TX})\oplus 
d\vartheta^2$ on $TY=T^HY\oplus TS^1$. 
We will denote by $\db_b^*$ the formal adjoint of $\db_b$ with 
respect to this metric and form the Kohn-Laplacian
\begin{equation}\label{kohn-lap}
\Box_b=\db_b\db_b^*+\db_b^*\db_b.
\end{equation}
The operators $\db_b^*$ and $\Box_b$ also commute with the action of $S^1$ 
on $Y$.
Consider the space $\cC^\infty(Y)_p$ of smooth functions $f$ on $Y$ which 
transform under the action $(y,\vartheta)\mapsto e^{i\vartheta}y$ of $S^1$ 
according to the law
\begin{equation}\label{equi}
f(e^{i\vartheta}y)=e^{ip\vartheta}f(y).
\end{equation}
This space of functions can be identified naturally with the space of smooth 
sections $\Omega^{0,0}(X,L^p)$.
More generally, the space of sections $\Omega^{0,k}(Y)_p$ 
which transform 
under the action of $S^1$ according to the law \eqref{equi} can be naturally 
identified with the space $\Omega^{0,k}(X,L^p)$.
Therefore, for each integer $p$, we get a subcomplex 
$\big(\Omega^{0,\sbullet}(Y)_p,\db_b\big)$ of the tangential Cauchy-Riemann 
complex \eqref{CRtangential}, isomorphic to the Dolbeault complex 
$\big(\Omega^{0,\sbullet}(X,L^p),\db^{L^p}\big)$.
Moreover, the action of $\Box_b$ on $\Omega^{0,\sbullet}(Y)_p$ is identical to 
the action of the Kodaira-Laplacian $\Box^{L^p}$ on
$\Omega^{0,\sbullet}(X,L^p)$, via the the complex isomorphism just mentioned.

Let us consider the spaces of $\Box_b$-harmonic spaces 
\begin{equation}\label{harmonic}
 \cH^{0,k}(Y)=\ker \Box_b|_{\Omega^{0,k}(Y)}\,,\quad 
 \cH^{0,k}(Y)_p=\ker \Box_b|_{\Omega^{0,k}(Y)_p}\,.
\end{equation}
Then
\begin{equation}\label{harmonic-desc}
 \cH^{0,k}(Y)=\oplus_{p\in\bZ} \cH^{0,k}(Y)_p\cong 
\oplus_{p\in\bZ} H^{0,k}(X,L^p)
\end{equation}
where $H^{0,k}(X,L^p):=\ker \Box^{L^p}|_{\Omega^{0,k}(X,L^p)}$ is the space 
of harmonic forms \eqref{f12}.

The {\em Szeg\"o projector\/} $\Pi^{0,k}$ is the orthogonal projection
from $\Omega^{0,k}(Y)$ to $\cH^{0,k}(Y)$.
The Szeg\"o projector $\Pi^{0,k}$ is finite dimensional for the degrees
$k\neq q,n-q$. This follows from the decomposition \eqref{harmonic-desc} and the 
vanishing theorem of Andreotti-Grauert \eqref{af12} applied for both $L$ and $L^*$.
It shows that there exits $p_0\in\bN$ such that $H^{0,k}(X,L^p)=0$ for all
$p\in\bZ$ with $|p|\geqslant p_0$ and all $k\neq q,n-q$ (note that $R^{L^*}=-R^L$
has $n-q$ negative and $q$ positive eigenvalues).

On the other hand, the Szeg\"o projector $\Pi^{0,k}$ is infinite dimensional 
for the degrees $k=q,n-q$ as shown by \eqref{af12} combined with 
the Riemann-Roch-Hirzebruch formula, which in turn is a consequence of 
the integration of the asymptotic expansion from
Theorem \ref{nonkahler} over the manifold $X$. To obtain the result for $k=n-q$
we have to replace $L$ by $L^*$ in the above mentioned results.

The description of the dimension of the harmonic spaces is consistent with 
the general geometric information from \cite{AnHi1,AnHi2} 
and \cite[p. 626]{Bou74}, where general hypersurfaces
$Y$ are considered, with non-degenerate Levi form of signature $(q,n-q)$.

The relation between the Bergman kernels $P^{0,q}_p$ considered in 
Theorem \ref{nonkahler} and the Szeg\"o kernel $\Pi^{0,q}$  is given by
\begin{equation}\label{Fourier}
P^{0,q}_p(x,x)=\frac{1}{2\pi}\int_{S^1}\Pi^{0,q}(e^{i\vartheta}y,y)
e^{-ip\vartheta}\,d\vartheta
\end{equation}
where $x\in X$ and $y\in Y$ satisfy $\pi(y)=x$. This means that
the $P^{0,q}_p(x,x)$ represent the Fourier coefficients of 
the distribution $\Pi^{0,q}(y,y)$.
Since we know the asymptotic expansion of $P^{0,q}_p(x,x)$ 
as $p\to\infty$ given in Theorem \ref{nonkahler}
we can recover from \eqref{Fourier} the restriction on the diagonal 
of the Szeg\"o kernel $\Pi^{0,q}$.

It could also be possible to work in the opposite direction and start with
 the Szeg\"o kernel. Namely, using a similar analysis as the one of 
Boutet de Monvel and Sj\"ostrand \cite{BouSj76} one can find the parametrix 
of the Szeg\"o kernel $\Pi^{0,q}$ and determine its singularity 
on the diagonal. Then, working as Zelditch \cite{Ze98} 
(where the case $q=0$ is considered) one can deduce the asymptotic
 of $P^{0,q}_p(x,x)$ for $p\to\infty$.

The same discussion applies to $\Pi^{0,n-q}$ and the Bergman kernels 
$P^{0,n-q}_p$ associated to $L^*$.

\section{The coefficient $\pmb{b}_1$} \label{s2}

This Section is organized as follows. In Section \ref{s2.1}, 
we state our main result, the formula for the coefficient $\pmb{b}_1$
for the spin$^c$ Dirac operator. 
In Section \ref{s2.2}, we obtain an asymptotic expansion of
the rescaled spin$^c$ Dirac operator $L^t_2$ (cf. \eqref{c27}) in normal coordinates.
In Section \ref{s2.3}, we finally compute the coefficient $\pmb{b}_1$.

\subsection{Main result}\label{s2.1}
We use the notation in Section \ref{s5.1}.
We denote by $I_{\bC\otimes E}$  the
projection from $\Lambda (T^{*(0,1)}X)\otimes E$ onto $\bC\otimes E$ 
under the decomposition
$\Lambda (T^{*(0,1)}X)= \bC \oplus \Lambda ^{>0} (T^{*(0,1)}X)$.
For any tensor $\psi$ on $X$, we denote by $\nabla ^{X}\psi$ 
the covariant derivative of $\psi$ induced by $\nabla ^{TX}$.
 Thus  $\nabla ^{X} J \in
T^*X \otimes \End(TX)$, $\nabla ^{X} \nabla ^{X} J\in T^*X
\otimes T^*X \otimes \End(TX)$.
Let $\{w_j\}$ be an orthonormal basis of $(T^{(1,0)}X, g^{TX})$, and
its dual basis $\{w ^j\}$. 
Let $\{e_i\}_{1\leqslant i\leqslant 2n}$ be an orthonormal basis of $(TX, g^{TX})$. 
Then we denote by $|\nabla ^X J |^2:= \sum_{ij}|(\nabla ^X_{e_i} J)e_j |^2$.

The following result is the main result of this paper.
\begin{thm} \label{t0.2} If $\bJ =J$, then for $\pmb{b}_1$ in \eqref{a0.3}, we have 
\begin{equation}\label{c1}
\begin{split}
\pmb{b}_1(x) &= \frac{1}{8\pi}\Big[r^X + \frac{1}{4} |\nabla ^X J |^2
+ 4 R^E (w_j,\ov{w}_j)\Big]I_{\bC\otimes E}
-\frac{1}{144\pi} \sum_{kl}|(\nabla ^X_{w_k} J)w_l |^2 I_{\bC\otimes E}\\
&+ \frac{1}{288\pi} \left \langle (\nabla ^X_{\ov{w}_k} J)\ov{w}_l,
\ov{w}_m \right \rangle 
\left \langle (\nabla ^X_{w_k} J)w_i,w_j \right \rangle
 \ov{w}^l \wedge \ov{w}^m I_{\bC\otimes E} i_{\ov{w}_j}\wedge i_{\ov{w}_i}\\
&- \frac{1}{8\pi } \left(\frac{1}{3} \left \langle R^{TX} w_i, \ov{w}_i\right \rangle
+ R^E\right) (\ov{w}_l,\ov{w}_m) \ov{w}^l\wedge \ov{w}^m I_{\bC\otimes E} \\
&+ \frac{1}{8\pi } \left(\frac{1}{3} \left \langle R^{TX} w_i, \ov{w}_i\right \rangle
+ R^E\right) (w_l,w_m) I_{\bC\otimes E} i_{\ov{w}_m}\wedge i_{\ov{w}_l}.
\end{split}
\end{equation}
Especially,  
\begin{align}\label{c2}
\tr|_{\Lambda (T^{*(0,1)}X)} [\pmb{b}_1(x)]
=  \frac{1}{8\pi}\Big[r^X + \frac{1}{4} |\nabla ^X J |^2
+ 4 R^E (w_j,\ov{w}_j)\Big].
\end{align}
\end{thm}
The term $r^X + \frac{1}{4} |\nabla ^X J |^2$ 
in \eqref{c2} is called the Hermitian scalar curvature in the literature
\cite{D97}, \cite[Chap.\,10]{Ga04},  \cite{MM04a} 
and is a natural substitute for the 
Riemannian scalar curvature in the almost-K\"ahler case. 
It was used by Donaldson \cite{D97} to define the moment map on the space of 
compatible almost-complex structures.

\subsection{Taylor expansion of the operator $L^t_2$}\label{s2.2}
To compare with \cite[\S 1.2]{MM04a}, in this part, we assume that 
$\om(\cdot,J\cdot)$ is positive, i.e. $q=0$ in \eqref{diag1}.

We fix $x_0\in X$. From now on, we identify $B^{T_{x_0}X}(0,\var)$
with $B^{X} (x_0,\var)$ by the exponential map 
$T_{x_0}X\ni Z\to \exp^X_{x_0}(Z)\in X$.
We identify $L_Z, E_Z$ and $(E_p)_Z$
for $Z\in B^{T_{x_0}X}(0,\var)$ to $L_{x_0}, E_{x_0}$ and $(E_p)_{x_0}$
by parallel transport with respect to the connections
$\nabla ^L, \nabla ^E$ and
 $\nabla^{E_p}$ along the curve $\gamma_Z :[0,1]\ni u \to uZ$.
 Let $\{e_i\}_i$ be an oriented orthonormal
basis of $T_{x_0}X$.
We also denote by $\{e ^i\}_i$ the dual basis of $\{e_i\}$.
 Let $\wi{e}_i (Z)$  be the parallel
transport of ${e}_i$ with respect to $\nabla^{TX}$
 along the above curve.

Let $S_L$ be an unit vector of $L_{x_0}$. 
Using $S_L$ and  the above discussion, we
get an isometry 
$E_{p}\simeq (\Lambda ( T^{*(0,1)}X)\otimes E)_{x_0}=: \bE_{x_0}$
on $B^{T_{x_0}X}(0,\var)$. Under our identification, 
$h^{E_p}$ is $h^{\bE_{x_0}}$ on $B^{T_{x_0}X}(0,\var)$.


For $s\in \cC^{\infty}( T_{x_0}X, \bE_{x_0}) $, set
\begin{align}\label{u0}
&\|s\|_{0,0}^2 = \int_{\bR^{2n}} 
|s(Z)|^2_{h^{\Lambda ( T^{*(0,1)}X)\otimes E}_{x_0}} dv_{TX}(Z).
\end{align}

Denote by  $\nabla_U$ the ordinary differentiation
 operator on $T_{x_0}X$ in the direction $U$.
If $\alpha = (\alpha_1,\cdots, \alpha_{2n})$ is a multi-index,
set $Z^\alpha = Z_1^{\alpha_1}\cdots Z_{2n}^{\alpha_{2n}}$.
Let $(\partial ^\alpha R^L)_{x_0}$
 be the tensor $(\partial^\alpha R^L)_{x_0}(e_i,e_j)
=\partial ^\alpha( R^L(e_i,e_j))_{x_0}$. 
We denote by
$\mR= \sum_i Z_i e_i =Z$ the radial vector field on $\bR^{2n}$.
Recall that the function $\kappa$ was defined in \eqref{c22}.
For  $s \in \cC^{\infty}(B^{T_{x_0}X}(0,\var), \bE_{x_0})$ 
and $Z\in B^{T_{x_0}X}(0,\var)$,
for $t=\frac{1}{\sqrt{p}}$, set
\begin{equation}\label{c27}
\begin{split}
&(S_{t} s ) (Z) =s (Z/t)\,,\quad    
\nabla_{t}=  S_t^{-1} t \kappa ^{\frac{1}{2}} 
\nabla ^{E_p}\kappa ^{-\frac{1}{2}} S_t\,, \\
&\nabla_{0,\,\pmb{\cdot}}=\nabla _{\pmb\cdot}+ \frac{1}{2}  R^L_{x_0}(\mR, \pmb\cdot)\,,
 \quad   L^t_2= S_t^{-1}  t^2 \kappa ^{\frac{1}{2}}
 D_p^{2}\kappa ^{-\frac{1}{2}} S_t\,.
\end{split}
\end{equation}
By our trivialization, $L^t_2$ is self-adjoint with respect to 
$\|\,\cdot\,\|_{0,0}$ on $\cC^{\infty}_0(B^{T_{x_0}X}(0,\var/t), \bE_{x_0})$.
Note that comparing with \cite[(4.37)]{DLM04a}, 
we conjugate with $\kappa^{1/2}$ in \eqref{c27}, 
which simplifies the computation of the coefficient $\pmb{b}_1$.

We  adopt
the convention that all tensors will be evaluated at the base point
$x_0\in X$, and most of the time, we will omit the subscript
$x_0$. Let  $\cL_0, \mO_1, \mO_2$ be the operators defined in 
\cite[Theorem 1.4]{MM04a} associated the renormalized Bochner-Laplacian
$\Delta_{p,0}$. Recall that $\tau= \sum_j R^L (w_j,\overline{w}_j)$,
$\om_d= - R^L (w_l,\overline{w}_m)\,\overline{w}^m\wedge
\,i_{\overline{w}_l}$.
Thus we have 
\begin{equation}\label{1c31}
\begin{split}
&\cL_0 = -\sum_j (\nabla_{0,e_j})^2
-\tau_{x_0},\\
&\mO_1(Z)= -\frac{2}{3}   ( \partial_j R^L)_{x_0} (\mR,e_i)Z_j
\nabla_{0,e_i}
 -\frac{1}{3} (\partial_i R^L)_{x_0} (\mR,
e_i) -(\partial_\mR \tau)_{x_0}\,,\\
&\mO_2(Z)=    \frac{1}{3} \left \langle R^{TX}_{x_0} (\mR,e_i)
\mR, e_j\right \rangle_{x_0}\nabla_{0,e_i}\nabla_{0,e_j}\\
 &\hspace*{3mm} +\Big [\frac{2}{3} \left \langle R^{TX}_{x_0}
(\mR, e_j) e_j,e_i\right \rangle 
- \Big(\frac{1}{2}\sum_{|\alpha|=2}(\partial ^{\alpha}R^L)_{x_0} 
\frac{Z^\alpha}{\alpha !} + R^E_{x_0} \Big)(\mR,e_i)\Big ]\nabla_{0,e_i}\\
&\hspace*{3mm} -\frac{1}{4} \nabla_{e_i}\Big(\sum_{|\alpha|=2}
(\partial ^{\alpha}R^L)_{x_0} \frac{Z^\alpha}{\alpha !}(\mR,e_i)\Big)
 -\frac{1}{9}\sum_i 
\Big[\sum_j (\partial_j R^L)_{x_0} (\mR,e_i)Z_j\Big]^2\\
&\hspace*{3mm}-\frac{1}{12}\Big[\cL_0, \left \langle R^{TX}_{x_0} (\mR,e_i)
\mR, e_i\right \rangle_{x_0} \Big] 
-\sum_{|\alpha|=2}(\partial ^{\alpha}\tau)_{x_0} 
\frac{Z^\alpha}{\alpha !} \,.
\end{split}
\end{equation}

\begin{thm}\label{t3.3} There are second order
 differential operators $L^0_2, \underline{\mQ}_r (r\geqslant 1)$ 
which are self-adjoint with respect to $\|\,\cdot\,\|_{0,0}$ 
on $\cC_0^\infty (\bR^{2n}, \bE_{x_0})$, and  
\begin{equation}\label{1c30}
\begin{split}
&L^0_2=\cL_0 - 2 \om_{d,x_0},\\
&\underline{\mQ}_1 = \mO_1
-\pi\sqrt{-1}\left \langle(\nabla_{\mR}^X \bJ)_{x_0}e_l, e_m\right \rangle
\,c(e_l)\,c(e_m)   +(\partial_\mR \tau)_{x_0}\,,\\
&\underline{\mQ}_2=  \mO_2 
- R^{\mathrm{Cliff}}_{x_0}(\mR,e_l)\nabla_{0,e_l}
-\frac{\pi}{2}\sqrt{-1}
\left \langle(\nabla^X\nabla^X \bJ)_{(\mR,\mR), x_0}e_l, e_m\right \rangle
\,c(e_l)\,c(e_m)\\
&\hspace{10mm}+  \frac{1}{2}  \Big( R^E_{x_0}
+\tfrac{1}{2}\tr\left[R^{T^{(1,0)}X}_{x_0}\right]\Big)(e_l,e_m)
\,c(e_l)\,c(e_m) 
+\sum_{|\alpha|=2}(\partial ^{\alpha}\tau)_{x_0} 
\frac{Z^\alpha}{\alpha !} + \frac{1}{4}r^X _{x_0}\,.
\end{split}
\end{equation}
 such that 
\begin{equation}\label{c30}
L^t_2= L^0_2 
 + \sum_{r=1}^\infty\underline{\mQ}_r t^r \,. 
\end{equation}
\end{thm}
\begin{proof} 
Set $g_{ij}(Z)= g^{TX}(e_i,e_j)(Z) =  \langle e_i,e_j\rangle_Z$ 
and let $(g^{ij}(Z))$ be the inverse of
the matrix $(g_{ij}(Z))$.
By \cite[Proposition 1.28]{BeGeVe} (cf. \cite[Lemma 4.5]{DLM04a}), we have
\begin{equation}\label{0c30}
\begin{split}
&g_{ij}(Z) =\delta_{ij} +  \frac{1}{3}
\left \langle R^{TX}_{x_0} (\mR,e_i) \mR, e_j\right \rangle_{x_0}
 + \cO (|Z|^3),\\
&\kappa(Z)= |\det (g_{ij}(Z))|^{1/2}  = 1 + 
\frac{1}{6} \left \langle R^{TX}_{x_0} (\mR,e_i) \mR, e_i\right \rangle_{x_0}
 + \cO (|Z|^3).
\end{split}
\end{equation}
If $\Gamma _{ij}^l$ is the connection form of $\nabla ^{TX}$
with respect to the basis $\{e_i\}$, we have $(\nabla
^{TX}_{e_i}e_j)(Z) = \Gamma _{ij}^l (Z) e_l$. 
Owing to \eqref{0c30} (cf. \cite[(1.32)]{MM04a}),
\begin{equation}\label{0c31}
\Gamma _{ij}^l (Z)= \frac{1}{3}
\left \langle R^{TX}_{x_0} (\mR, e_j) e_i
+  R^{TX}_{x_0} (\mR, e_i) e_j, e_l\right \rangle_{x_0}
 + \cO(|Z|^2). 
\end{equation}

Let $\Gamma ^E$, $\Gamma ^L$,  $\Gamma ^{\text{Cliff}}$
be the connection forms of $\nabla^E$, $\nabla^L$, $\nabla^{\text{Cliff}}$
 with respect to any fixed frames for $E$, $L$, $\Lambda (T^{*(0,1)}X)$
 which are parallel along the curve $\gamma_Z$
under our trivializations on $B^{T_{x_0}X}(0,\var)$.
By \cite[Proposition 1.18]{BeGeVe}
the Taylor coefficients of $\Gamma^{\sbullet} (e_j) (Z)$ at $x_0$
to order $r$ are only determined by those of $R^{\sbullet}$ to order $r-1$, and
\begin{align}\label{0c39}
\sum_{|\alpha|=r}  (\partial^\alpha
 \Gamma^{\sbullet} ) _{x_0} (e_j) \frac{Z^\alpha}{\alpha !}
=\frac{1}{r+1} \sum_{|\alpha|=r-1}
(\partial^\alpha R^{\sbullet} ) _{x_0}(\mR, e_j) 
  \frac{Z^\alpha}{\alpha !}.
\end{align}
\eqref{c27}, \eqref{0c39} yield on $B^{T_{x_0}X}(0, \var/t)$,
\begin{equation}\label{0c36}
\begin{split}
\nabla_{t, e_i}|_{Z} &= \kappa ^{\frac{1}{2}}(tZ)\Big[\nabla_{e_i}
+ \Big(\frac{1}{t} \Gamma ^L (e_i) + t \Gamma ^E (e_i)
+ t \Gamma ^{\text{Cliff}}(e_i)\Big)(tZ) \Big]
\kappa ^{-\frac{1}{2}}(tZ)\\
&=\kappa ^{\frac{1}{2}}(tZ)\Big[\nabla_{e_i}
+\Big( \frac{1}{2} R^L_{x_0} 
+ \frac{t}{3} (\partial_k R^L)_{x_0} Z_k \\
&\hspace*{10mm} +\frac{t^2}{4}\sum_{|\alpha|=2}
(\partial ^{\alpha}R^L)_{x_0} \frac{Z^\alpha}{\alpha !}
+ \frac{t^2}{2} R^E_{x_0}
+ \frac{t^2}{2} R^{\text{Cliff}}_{x_0} \Big )  (\mR,e_i)+ \cO (t^3)
\Big] \kappa ^{-\frac{1}{2}}(tZ) . 
\end{split}
\end{equation}

 By the definition of $\nabla^{{\rm Cliff}}$, 
for $Y,U\in \cC^\infty(X, TX)$,
\be\label{c31}
[\nabla^{\text{Cliff}}_U, c(Y)] = c(\nabla ^{TX}_U Y).
\ee
For $\psi\in T^*X \otimes \End(\Lambda (T^{*(0,1)}X))
\simeq T^*X \otimes (C(TX)\otimes_\bR \bC)$, where $C(TX)$ is the Clifford
algebra bundle of $TX$, we still denote by $\nabla^X\psi$ the covariant 
derivative of $\psi$ induced by $\nabla^{TX}$.
 By using \eqref{c31}, we observe that
\begin{align}\label{0c38}
\nabla^X _Y (\psi(\wi{e}_j)c(\wi{e}_j))
&= (\nabla^X _Y \psi)(\wi{e}_j)c(\wi{e}_j) 
+ \psi(\nabla^{TX}_Y\wi{e}_j)c(\wi{e}_j) 
+ \psi(\wi{e}_j)c(\nabla^{TX}_Y\wi{e}_j))\\
&= (\nabla^X _Y \psi)(\wi{e}_j)c(\wi{e}_j) .\nonumber
\end{align}
Thus for $k\geqslant 2$,
\begin{equation}\label{0c38a}
\begin{split}
&\left(R^L(\wi{e}_l,\wi{e}_m) \,c(\wi{e}_l)\,c(\wi{e}_m)\right)(tZ)=\\
&= \sum_{r=0}^k \frac{\partial^r}{\partial t^r}
\left[(R^L(\wi{e}_l,\wi{e}_m) \,c(\wi{e}_l)\,c(\wi{e}_m))(tZ)\right] |_{t=0}\, \frac{t^r}{r !} + \cO(t^{k+1})\\
&=\Big (R^L_{x_0}+ t (\nabla^X _\mR R^L)_{x_0} 
+ \frac{t^2}{2} (\nabla^X \nabla^X  R^L)_{(\mR,\mR),x_0} \Big ) 
(e_l,e_m)  \,c(e_l)\,c(e_m) + \cO(t^{3}).
\end{split}
\end{equation}
 
Owing to  \eqref{Lich}, \eqref{c27}, \eqref{0c36}
 and \eqref{0c38}-\eqref{0c38a}, 
\begin{equation}\label{0c41}
\begin{split}
 L^t_2 =& - g^{ij}(tZ) \left(\nabla_{t, e_i}\nabla_{t, e_j}
- t\, \Gamma ^l_{ij}(tZ) \nabla_{t, e_l} \right )+  \frac{t^2}{4} r^X (tZ)\\
&+ \frac{1}{2} \left\{ \Big[ R^L
+  t^2 \left( R^E +\tfrac{1}{2}\tr [R^{T^{(1,0)}X} ]\right)
\Big] (\wi{e}_l,\wi{e}_m)  \,c(\wi{e}_l)\,c(\wi{e}_m)\right\} (tZ) \\
&= - g^{ij}(tZ) \left(\nabla_{t, e_i}\nabla_{t, e_j}
- t \Gamma ^l_{ij}(tZ) \nabla_{t, e_l} \right )  \\
 &+\frac{1}{2} \Big( R^L_{x_0} + t (\nabla^X_{\mR} R^L) _{x_0}
+\frac{t^2}{2} (\nabla^X \nabla^X  R^L)_{(\mR,\mR),x_0} \Big ) 
 (e_l,e_m)\,c(e_l)\,c(e_m)\\
&+ \frac{t^2}{2}\left( R^E_{x_0} 
+\tfrac{1}{2}\tr\left[R^{T^{(1,0)}X}_{x_0} \right]\right)
(e_l,e_m) \,c(e_l)\,c(e_m)
+ \frac{t^2}{4} r^X _{x_0} 
+ \cO (t^3).
\end{split}
\end{equation}

Note that $\nabla^X g^{TX}=0$. 
Comparing \eqref{a0.1}, \eqref{0c30}, \eqref{0c36}, \eqref{0c41},  
with \cite[(1.37)]{MM04a}, we get \eqref{c30}.

To prove the self-adjointness of $L^0_2,$ and $\underline{\mQ}_r (r\geqslant 1)$
with respect to $\|\,\cdot\,\|_{0,0}$ on $\cC_0^\infty (\bR^{2n}, \bE_{x_0})$, we
observe that it follows from the fact that $L^t_2$ is self-adjoint with respect to $\|\,\cdot\,\|_{0,0}$ 
on $\cC^{\infty}_0(B^{T_{x_0}X}(0,\var/t), \bE_{x_0})$.
\end{proof}

 Let $P^N$, $\underline{P}^N$  be the orthogonal projections from  
$(L^2 (\bR^{2n}, \bE_{x_0}), \|\quad\|_{0,0})$
 onto $N=\Ker \cL_0$, $\Ker L^0_2$,
and let $P^N(Z,Z')$, $\underline{P}^N(Z,Z')$
 be the smooth kernel of $P^N$, $\underline{P}^N$ with respect 
to $dv_{TX}(Z)$. Set $P^ {N^\bot}= \Id -P^N$, 
$\underline{P}^{N^\bot}= \Id -\underline{P}^N$.
Recall that $q=0$ in \eqref{diag1}, 
thus $\om_d \leq -\mu_0$ on $\Lambda^{>0}(T^{*(0,1)}X)$, by \eqref{1c30},
\begin{align}\label{0ue62}
\underline{P}^N(Z,Z') = P^N(Z,Z') I_{\bC\otimes E}.
\end{align}

\begin{thm}\label{t3.5} We have the relation
\begin{align}\label{0c43}
\underline{P}^N \underline{Q}_1 \underline{P}^N=0.
\end{align}
\end{thm}
\begin{proof} By \cite[(1.92), (1.94), (1.96)]{MM04a}, we have 
\begin{equation}\label{0c44}
\begin{split}
& (\nabla ^{X}_U R^L) (V,W)
= -2 \pi \sqrt{-1}\langle (\nabla ^{X}_U\bJ) V,W \rangle,\\
&(\nabla _\mR\tau)_{x_0} =-2\pi \sqrt{-1}
\left\langle(\nabla ^{X}_{\mR}\bJ) w_i, \ov{w}_i\right\rangle,\\
&P^N\mO_1P^N=0. 
\end{split}
\end{equation}
By \eqref{1c30}, \eqref{0ue62} and \eqref{0c44}, we get \eqref{0c43}.
\end{proof}

Set 
\begin{equation}\label{g51}
\begin{split}
 F_{2}=&(L^0_2)^{-1}\underline{P}^{N^\bot}\underline{\mQ}_1 (L^0_2)^{-1}\underline{P}^{N^\bot}\underline{\mQ}_1 \underline{P}^N
- (L^0_2)^{-1}\underline{P}^{N^\bot}\underline{\mQ}_2 \underline{P}^N\\
&+ \underline{P}^N\underline{\mQ}_1(L^0_2)^{-1}\underline{P}^{N^\bot}\underline{\mQ}_1(L^0_2)^{-1}\underline{P}^{N^\bot}
- \underline{P}^N \underline{\mQ}_2(L^0_2)^{-1}\underline{P}^{N^\bot}\\
&+ \underline{P}^{N^\bot}(L^0_2)^{-1}\underline{\mQ}_1 \underline{P}^N
\underline{\mQ}_1 (L^0_2)^{-1}\underline{P}^{N^\bot}
-\underline{P}^N\underline{\mQ}_1\underline{P}^{N^\bot}
(L^0_2)^{-2}\underline{\mQ}_1 \underline{P}^N.
\end{split}
\end{equation}
Then by Theorem \ref{t3.5} and
 the same argument as in \cite[\S 1.5, 1.6]{MM04a}, we get
\begin{equation}\label{g52}
\pmb{b}_1(x_0)= F_{2}(0,0). 
\end{equation}
By Theorem \ref{t3.3}, the third and fourth terms in \eqref{g51}
are adjoint of the first two terms, thus we only need to compute 
the first two terms and the last two terms in \eqref{g51}.

\subsection{Computing the coefficient $\pmb{b}_1$} \label{s2.3}
From now on, we assume that $\bJ=J$. By \cite[(2.13)]{MM04a},
\begin{equation}\label{g29}
\begin{split}
&\text{ $\nabla ^{X}_U J$ is skew-adjoint
and the tensor $\left\langle(\nabla ^{X}_\cdot J)\cdot,
\cdot \right\rangle$  is of the type}\\ 
&(T^{*(1,0)}X)^{\otimes 3} \oplus (T^{*(0,1)}X)^{\otimes 3}.
\end{split}
\end{equation}

In what follows we will use the complex coordinates $z=(z_1,\cdots,z_n)$,
such that $w_i=\sqrt{2}\tfrac{\partial}{\partial z_i}$
is an orthonormal basis of  $T^{(1,0)}_{x_0} X$.
Then $Z=z+\overline{z}$ and 
 will also identify $z$ to $\sum_i z_i\tfrac{\partial}{\partial z_i}$ 
and $\overline{z}$ to
$\sum_i\overline{z}_i\tfrac{\partial}{\partial\overline{z}_i}$ when
we consider $z$ and $\overline{z}$ as vector fields. Remark that
\begin{equation}\label{g0}
\Big\lvert\tfrac{\partial}{\partial z_i}\Big\rvert^2=
\Big\lvert\tfrac{\partial}{\partial\overline{z}_i}\Big\rvert^2
=\dfrac{1}{2}\,,
\quad\text{so that $|z|^2=|\overline{z}|^2=\dfrac{1}{2} |Z|^2$\,.}
\end{equation}
It is very useful to rewrite $\cL_0$ by using  
the creation and annihilation operators. Set
\begin{equation}\label{0g1}
\quad b_i=-2\nabla_{0,\tfrac{\partial}{\partial z_i}}
=-2{\tfrac{\partial}{\partial z_i}}+\pi \overline{z}_i,\quad
b^{+}_i=2\nabla_{0,\tfrac{\partial}{\partial \overline{z}_i}}
=2{\tfrac{\partial}{\partial\overline{z}_i}}+\pi  z_i,
\quad b=(b_1,\cdots,b_n)\,.
\end{equation}
Then for any polynomial  $g(z,\ov{z})$ on $z$ and $\ov{z}$, 
\begin{equation}\label{g2}
\begin{split}
&[b_i,b^{+}_j]=b_i b^{+}_j-b^{+}_j b_i =-4 \pi  \delta_{i\,j},\\
&[b_i,b_j]=[b^{+}_i,b^{+}_j]=0\, ,\\
& [g(z,\ov{z}),b_j]=  2 \tfrac{\partial}{\partial z_j}g(z,\ov{z}), 
\quad  [g(z,\ov{z}),b_j^+]
= - 2\tfrac{\partial}{\partial \overline{z}_j}g(z,\ov{z})\,. 
\end{split}
\end{equation}
As $\bJ=J$, $a_j=2\pi$ in \eqref{diag2}.
 By \eqref{1c31} and \eqref{1c30} (cf. \cite[(2.23)]{MM04a}), 
\begin{equation}\label{g3}
\begin{split}
&\mO_1=-\frac{4\pi\sqrt{-1}}{3}b_i 
\left\langle(\nabla ^{X}_{\overline z} J)\overline z,
\tfrac{\partial}{\partial\ov{z}_i}\right\rangle
+ \frac{4\pi\sqrt{-1}}{3} \left\langle(\nabla ^{X}_z J)z,
\tfrac{\partial}{\partial z_i} \right\rangle b_i^+,\\
&\cL_0= b_j b^{+}_j,\quad  
L^0_2=  b_j b^{+}_j 
+ 4\pi   \ov{w}^j\wedge i_{\ov{w}_j}. 
\end{split}
\end{equation}
We found the following result in \cite[Theorem 1.15]{MM04a}:
\begin{thm}\label{t3.4}
The spectrum of the restriction of $\cL_0$ on $L^2(\bR^{2n})$ is given by
\begin{equation}\label{0g4}
\spec{{\cL_0}|_{L^2(\bR^{2n})}}=
\Big\lbrace  4\pi \sum_{i=1}^n\alpha_i \,:\, 
 \alpha =(\alpha_1,\cdots,\alpha_n)\in\bN ^n\Big\rbrace
\end{equation}
and an orthogonal basis of the eigenspace of $4\pi \sum_{i=1}^n\alpha_i$
is given by
\begin{equation}\label{0g5}
b^{\alpha}\big(z^{\beta}\exp\big({-    \frac{\pi}{2}\sum_i
|z_i|^2}\big)\big)\,,\quad\text{with $\beta\in\bN^n$}\,.
\end{equation}
\end{thm}
From \eqref{0g5}, we get
\begin{equation}\label{0g6}
\begin{split}
P^N(Z,Z') &= \exp\Big(-\frac{\pi}{2}\sum_i
\big(|z_i|^2+|z^{\prime}_i|^2 -2z_i\overline{z}_i'\big)\Big).
\end{split}
\end{equation}

As ${\bf J}=J$ here, we know the function $\tau$ therein is $2\pi n$.
By \eqref{1c30}, \eqref{g29}-\eqref{0g1}, 
\begin{equation}\label{g4}
\begin{split}
\underline{\mQ}_1=&\mO_1 
-2 \pi \sqrt{-1}  \Big [ 
 \left\langle (\nabla ^{X}_{\ov{z}}J) 
\tfrac{\partial}{\partial \ov{z}_l}, \tfrac{\partial}{\partial \ov{z}_m}
\right\rangle d\ov{z}_l d \ov{z}_m + 4  \left\langle (\nabla ^{X}_{z}J) 
\tfrac{\partial}{\partial z_l}, \tfrac{\partial}{\partial z_m}
\right\rangle i_{\tfrac{\partial}{\partial \ov{z}_l}}
i_{\tfrac{\partial}{\partial \ov{z}_m}} \Big], \\
\underline{\mQ}_2=&\mO_2
+ R^{\text{Cliff}}_{x_0}(\mR, \tfrac{\partial}{\partial \ov{z}_i}) b_i
- R^{\text{Cliff}}_{x_0}(\mR, \tfrac{\partial}{\partial z_i}) b_i^+\\
&\hspace{5mm}- \frac{\pi}{2} \sqrt{-1}
\left \langle (\nabla ^{X} \nabla ^{X}J)_{(\mR,\mR)}e_l,e_m\right \rangle
 \,c(e_l)\,c(e_m)\\
&\hspace{5mm}+\frac{1}{2} 
\Big( R^E_{x_0}+\frac{1}{2} \tr[R^{T^{(1,0)}X}_{x_0}]\Big)
(e_l,e_m) \,c(e_l)\,c(e_m)+ \frac{1}{4}r^X_{x_0}.
\end{split}
\end{equation}

Recall that by \cite[(1.98), (2.24)]{MM04a}, \eqref{0g1},
\begin{align}\label{g5}
&b^+_i P^N=0, \quad 
(b_i P^N)(Z,Z^\prime)= 2\pi (\ov{z}_i-\ov{z}_i^\prime) P^N(Z,Z^\prime), 
\quad (\mO_1  P^N)(Z,0)=0.
\end{align}
By \eqref{0ue62}, \eqref{g29}, \eqref{g3}, \eqref{0g6}, 
\eqref{g4} and \eqref{g5}, as in \cite[(2.24)]{MM04a},
\begin{multline}\label{g6}
(\underline{\mQ}_1 \underline{P}^N)(Z,Z^\prime)
=\left[ -\frac{2\sqrt{-1}}{3} b_ib_j \Big\langle 
(\nabla ^{X}_{\tfrac{\partial}{\partial\ov{z}_j}} J)\ov{z}^\prime,
\tfrac{\partial}{\partial\ov{z}_i}\Big\rangle 
-\frac{4\pi\sqrt{-1}}{3}  b_i
 \left\langle(\nabla ^{X}_{\ov{z}^\prime} J)\ov{z}^\prime,
\tfrac{\partial}{\partial\ov{z}_i}\right\rangle \right.\\
\left. - \sqrt{-1} \left\langle 
(\nabla ^{X}_{\tfrac{\partial}{\partial \ov{z}_k}}J) 
\tfrac{\partial}{\partial \ov{z}_l}, \tfrac{\partial}{\partial \ov{z}_m}
\right\rangle d\ov{z}_l d \ov{z}_m  (b_k + 2 \pi \ov{z}_k^\prime)
\right] P^N(Z,Z^\prime) I_{\bC\otimes E}.
\end{multline}
Thus by Theorem \ref{t3.4} and relations \eqref{g29},
 \eqref{g3} and \eqref{g6}, 
\begin{multline}\label{0g7}
((L^0_2)^{-1}\underline{P}^{N^\bot}
\underline{\mQ}_1 \underline{P}^N)(Z,Z^\prime)
=- \sqrt{-1} \left[ \frac{b_ib_j}{12\pi} \Big\langle 
(\nabla ^{X}_{\tfrac{\partial}{\partial\ov{z}_j}} J)\ov{z}^\prime,
\tfrac{\partial}{\partial\ov{z}_i}\Big\rangle 
+ \frac{b_i}{3} \left\langle(\nabla ^{X}_{\ov{z}^\prime} J)\ov{z}^\prime,
\tfrac{\partial}{\partial\ov{z}_i}\right\rangle \right.\\
\left. +
\left\langle (\nabla ^{X}_{\tfrac{\partial}{\partial \ov{z}_k}}J) 
\tfrac{\partial}{\partial \ov{z}_l}, \tfrac{\partial}{\partial \ov{z}_m}
\right\rangle d\ov{z}_l d \ov{z}_m  
\Big(\frac{b_k}{12\pi} +\frac{\ov{z}_k^\prime}{4}\Big)\right] 
 P^N(Z,Z^\prime) I_{\bC\otimes E}.
\end{multline}

By \eqref{g29},  \eqref{g5} and \eqref{0g7},
\begin{align}\label{g7}
\begin{split}
&(\underline{P}^{N^\bot}(L^0_2)^{-1}\underline{\mQ}_1
 \underline{P}^N)(0,Z^\prime) 
= - \frac{\sqrt{-1}}{12} 
\left\langle (\nabla ^{X}_{\ov{z}^\prime}J) 
\tfrac{\partial}{\partial \ov{z}_l}, \tfrac{\partial}{\partial \ov{z}_m}
\right\rangle  d\ov{z}_l d \ov{z}_m  P^N(0,Z^\prime) 
I_{\bC\otimes E},\\
&(\underline{P}^{N^\bot}(L^0_2)^{-1}\underline{\mQ}_1 \underline{P}^N)(Z,0) 
= - \frac{\sqrt{-1}}{6} 
\left\langle (\nabla ^{X}_{\ov{z}}J) 
\tfrac{\partial}{\partial \ov{z}_l}, \tfrac{\partial}{\partial \ov{z}_m}
\right\rangle  d\ov{z}_l d \ov{z}_m   P^N(Z,0) I_{\bC\otimes E}.
\end{split}
\end{align}
Recall that $L^0_2, \mQ_1$ are self-adjoint with respect 
to $\|\,\cdot\,\|_{0,0}$ on $\cC^\infty_0(\bR^{2n}, \bE_{x_0})$.
After taken the adjoint of \eqref{g7}, by \eqref{g0}, we get 
\begin{align}\label{0g8}
\begin{split}
&(\underline{P}^N\underline{\mQ}_1(L^0_2)^{-1}\underline{P}^{N^\bot})
(Z^\prime,0) 
= \frac{\sqrt{-1}}{3} 
\left\langle (\nabla ^{X}_{z^\prime}J)
\tfrac{\partial}{\partial z_l}, \tfrac{\partial}{\partial z_m}
\right\rangle  I_{\bC\otimes E} 
i_{\tfrac{\partial}{\partial \ov{z}_m}}i_{\tfrac{\partial}{\partial \ov{z}_l}}
P^N(Z^\prime,0) ,\\
&(\underline{P}^N\underline{\mQ}_1(L^0_2)^{-1}\underline{P}^{N^\bot}) (0,Z)
=  \frac{2\sqrt{-1}}{3\pi} 
\left\langle (\nabla ^{X}_{z}J)
 \tfrac{\partial}{\partial z_l}, \tfrac{\partial}{\partial z_m}
\right\rangle  I_{\bC\otimes E} 
i_{\tfrac{\partial}{\partial \ov{z}_m}}i_{\tfrac{\partial}{\partial \ov{z}_l}}
P^N(0, Z). 
\end{split}
\end{align}

Note that $\int_{\bC} |z|^2 e^{-\pi |z|^2} =1/\pi$.
By \eqref{g29}, \eqref{g0}, \eqref{0g6}, \eqref{g7} and \eqref{0g8},
\begin{align}\label{0g9}
&\Big(\underline{P}^{N^\bot}(L^0_2)^{-1}\underline{\mQ}_1 \underline{P}^N
\underline{\mQ}_1(L^0_2)^{-1}\underline{P}^{N^\bot}\Big)(0,0)\\
&\hspace*{10mm}= \frac{1}{36\pi}
\left\langle (\nabla ^{X}_{\tfrac{\partial}{\partial \ov{z}_k}}J) 
\tfrac{\partial}{\partial \ov{z}_l}, \tfrac{\partial}{\partial \ov{z}_m}
\right\rangle 
 \left\langle (\nabla ^{X}_{\tfrac{\partial}{\partial z_k}}J)
 \tfrac{\partial}{\partial z_i}, \tfrac{\partial}{\partial z_j}\right\rangle
 d\ov{z}_l d \ov{z}_m  I_{\bC\otimes E} 
i_{\tfrac{\partial}{\partial \ov{z}_j}}
i_{\tfrac{\partial}{\partial \ov{z}_i}},\nonumber\\
&- \Big(\underline{P}^N\underline{\mQ}_1\underline{P}^{N^\bot}
(L^0_2)^{-2}\underline{\mQ}_1 \underline{P}^N\Big) (0,0)\label{0g9a}\\
&\hspace*{7mm}=-  \frac{1}{9\pi} 
 \left\langle (\nabla ^{X}_{\tfrac{\partial}{\partial z_k}}J)
 \tfrac{\partial}{\partial z_i}, \tfrac{\partial}{\partial z_j}\right\rangle
 \left\langle (\nabla ^{X}_{\tfrac{\partial}{\partial \ov{z}_k}}J) 
\tfrac{\partial}{\partial \ov{z}_l}, \tfrac{\partial}{\partial \ov{z}_m}
\right\rangle 
I_{\bC\otimes E} i_{\tfrac{\partial}{\partial \ov{z}_j}}
i_{\tfrac{\partial}{\partial \ov{z}_i}}  
d\ov{z}_l d \ov{z}_m 
I_{\bC\otimes E} \nonumber\\
&\hspace*{7mm}= - \frac{1}{9\pi}\left\langle 
(\nabla ^{X}_{\tfrac{\partial}{\partial \ov{z}_k}}J) 
\tfrac{\partial}{\partial \ov{z}_l},
(\nabla ^{X}_{\tfrac{\partial}{\partial z_k}}J)
 \tfrac{\partial}{\partial z_l}\right\rangle . \nonumber
\end{align} 
Let $h_i(Z)$ (resp. $F(Z)$) be polynomials in $Z$ with degree $1$
(resp. $2$), then by Theorem \ref{t3.4}, 
\eqref{g2}, \eqref{g3} and \eqref{g5},
\begin{equation}\label{0g10}
\begin{split}
&(\cL_0^{-1}P^{N^\bot} h_i b_i P^N)(0,0) 
= (\cL_0^{-1} P^{N^\bot}b_i h_i P^N)(0,0)
= -\frac{1}{2\pi}\frac{\partial h_i}{\partial z_i},\\
&(\cL_0^{-1}P^{N^\bot} FP^N)(0,0) = -\frac{1}{4\pi ^2}
\frac{\partial ^2 F}{\partial z_i\partial \ov{z}_i},  \\
&((L^0_2)^{-1} F d\ov{z}_l d \ov{z}_m \underline{P}^N)(0,0)
= \frac{1}{2\pi} \Big((L^0_2)^{-1} 
\frac{\partial ^2 F}{\partial z_i\partial \ov{z}_i} (b_i z_i +2) 
 d\ov{z}_l d \ov{z}_m \underline{P}^N\Big)(0,0)\\
&\hspace*{7mm}
= \frac{1}{8\pi ^2} \frac{\partial ^2 F}{\partial z_i\partial \ov{z}_i}
\Big( (\frac{b_i z_i}{3} +1)  d\ov{z}_l d \ov{z}_m \underline{P}^N\Big)(0,0)
= \frac{1}{24\pi ^2}\frac{\partial ^2 F}{\partial z_i\partial \ov{z}_i} 
d\ov{z}_l d \ov{z}_m  I_{\bC\otimes E}.
\end{split}
\end{equation}
Observe that for a monomial $Q$ in $b_i$, $z_i$, $b_i^+, \ov{z}_i$,
if the total degree of $b_i$, $z_i$ is not the same as the total degree of 
 $b_i^+, \ov{z}_i$, then $(Q\underline{P}^{N})(0,0)=0$.
Theorem \ref{t3.4} and relations \eqref{g0}, \eqref{g3}, \eqref{g4}, 
\eqref{g5}, \eqref{g7} and \eqref{0g10} imply that
\begin{equation}\label{g8}
\begin{split}
((L^0_2)^{-1}\underline{P}^{N^\bot}&\underline{\mQ}_1
(L^0_2)^{-1}\underline{\mQ}_1  \underline{P}^N)(0,0)
=\Big( (L^0_2)^{-1}\underline{P}^{N^\bot}
(-8 \pi \sqrt{-1}) \left\langle (\nabla ^{X}_{z}J) 
\tfrac{\partial}{\partial z_i}, \tfrac{\partial}{\partial z_j}
\right\rangle\\
&i_{\tfrac{\partial}{\partial \ov{z}_i}}
i_{\tfrac{\partial}{\partial \ov{z}_j}}
\times   \frac{- \sqrt{-1}}{6} 
\left\langle (\nabla ^{X}_{\ov{z}}J) 
\tfrac{\partial}{\partial \ov{z}_l}, \tfrac{\partial}{\partial \ov{z}_m}
\right\rangle  d\ov{z}_l d \ov{z}_m   P^N\Big)(0,0)  I_{\bC\otimes E}\\
=& \frac{4\pi}{3}\Big( (L^0_2)^{-1}\underline{P}^{N^\bot}
\left\langle (\nabla ^{X}_{z}J) \tfrac{\partial}{\partial z_l}, 
 (\nabla ^{X}_{\ov{z}}J) 
\tfrac{\partial}{\partial \ov{z}_l}\right\rangle  P^N\Big)(0,0)  
I_{\bC\otimes E}\\
=& -\frac{1}{3\pi} \left\langle (\nabla ^{X}_{\tfrac{\partial}{\partial z_k}}J) 
\tfrac{\partial}{\partial z_l},  
(\nabla ^{X}_{\tfrac{\partial}{\partial \ov{z}_k}}J) 
\tfrac{\partial}{\partial \ov{z}_l}\right\rangle 
I_{\bC\otimes E}.
\end{split}
\end{equation}

Now we will compute 
$((L^0_2)^{-1}\underline{P}^{N^\bot}\underline{\mQ}_2\underline{P}^N)(0,0)$.
By \eqref{c20} and \eqref{g29}, 
\begin{align}\label{g13}
\tr [R^{T^{(1,0)}X}] = 
2 \left\langle R^{TX}\tfrac{\partial}{\partial z_j}, 
\tfrac{\partial}{\partial \ov{z}_j} \right\rangle
+ \frac{1}{2} \left\langle 
(\nabla ^X_{\tfrac{\partial}{\partial z_l}} J) \tfrac{\partial}{\partial z_j},
 (\nabla ^X_{\tfrac{\partial}{\partial \ov{z}_m}} J)
\tfrac{\partial}{\partial \ov{z}_j} \right\rangle d\ov{z}_m \wedge dz_l.
\end{align}

By definition, for $U,V\in TX$, we have
\begin{align}\label{g9}
& (\nabla ^{X}\nabla ^{X}J)_{(U,V)}- (\nabla ^{X}\nabla ^{X}J)_{(V,U)}
=[R^{TX}(U,V),J].
\end{align}
By \cite[(2.19)]{MM04a}, 
for $u_1,u_2, u_{3}\in T^{(1,0)}X$, $\ov{v}_1,\ov{v}_2\in T^{(0,1)}X$,
\begin{equation}\label{g35} 
\begin{split}
&(\nabla ^{X}\nabla ^{X}J)_{(u_1,u_2)}u_{3},\,
(\nabla ^{X}\nabla ^{X}J)_{(\ov{u}_1,\ov{u}_2)}u_{3}\in T^{(0,1)}X,
\quad   (\nabla ^{X}\nabla ^{X}J)_{(\ov{u}_1,u_2)}\ov{u}_{3} \in T^{(0,1)}X,\\
&2\sqrt{-1} \left\langle (\nabla ^{X}\nabla ^{X}J)_{(u_1,\ov{v}_1)}
u_2 ,  \ov{v}_2\right\rangle
= \left\langle  (\nabla ^{X}_{u_1}J) u_2,  
(\nabla ^{X}_{\ov{v}_1}J)\ov{v}_2\right\rangle \, .  
\end{split}
\end{equation}

By \eqref{c20}, \eqref{0ue62}, \eqref{0g6}, 
\eqref{g4}, \eqref{g5} and \eqref{g13},
\begin{equation}\label{g10}
\begin{split}
&\underline{P}^{N^\bot}(\underline{\mQ}_2- \mO_2)\underline{P}^N= 
\underline{P}^{N^\bot}\Big \{  \frac{1}{2}\tr\left[R^{T^{(1,0)}X}\right]
(\mR, \tfrac{\partial}{\partial \ov{z}_i}) b_i \\
&\hspace{15mm}-\left\langle \Big(R^{TX}(\mR, \tfrac{\partial}{\partial \ov{z}_i}) b_i 
-2\pi \sqrt{-1}  (\nabla ^{X} \nabla ^{X}J)_{(\mR,\mR)}\Big)
 \tfrac{\partial}{\partial z_l}, 
\tfrac{\partial}{\partial \ov{z}_l} \right\rangle  \\
&\hspace{15mm}+  \left\langle \Big(\frac{1}{2}R^{TX}(\mR, 
\tfrac{\partial}{\partial \ov{z}_i}) b_i 
-\pi \sqrt{-1}  (\nabla ^{X} \nabla ^{X}J)_{(\mR,\mR)}\Big) 
\tfrac{\partial}{\partial \ov{z}_l}, 
\tfrac{\partial}{\partial \ov{z}_m} \right\rangle 
d\ov{z}_l d\ov{z}_m\\
&\hspace{15mm}+\Big ( R^E+ \frac{1}{2} \tr\left[R^{T^{(1,0)}X}\right]\Big )
(\tfrac{\partial}{\partial \ov{z}_l}, 
\tfrac{\partial}{\partial \ov{z}_m}) 
d\ov{z}_l d\ov{z}_m \Big\}\underline{P}^N.
\end{split}
\end{equation}
Thus by \eqref{g5}, \eqref{0g10}, \eqref{g9} and \eqref{g10}, 
\begin{equation}\label{g11}
\begin{split}
-&\left((L^0_2)^{-1}\underline{P}^{N^\bot}(\underline{\mQ}_2- \mO_2)
\underline{P}^N\right)(0,0)=\left\{\frac{1}{4\pi} \tr\left[R^{T^{(1,0)}X}\right]
(\tfrac{\partial}{\partial z_i}, \tfrac{\partial}{\partial \ov{z}_i})\right.\\
 &- \frac{1}{2 \pi} \left\langle R^{TX}
(\tfrac{\partial}{\partial z_i}, \tfrac{\partial}{\partial \ov{z}_i})
-\sqrt{-1} \left( 2(\nabla ^{X} \nabla ^{X}J)_{(\tfrac{\partial}{\partial z_i},
 \tfrac{\partial}{\partial \ov{z}_i})}  
- [ R^{TX}(\tfrac{\partial}{\partial z_i}, 
\tfrac{\partial}{\partial \ov{z}_i}), J] \right)
\tfrac{\partial}{\partial z_l}, \tfrac{\partial}{\partial \ov{z}_l}
\right \rangle \\
&- \frac{1}{24\pi}\left\langle R^{TX}
(\tfrac{\partial}{\partial z_i}, \tfrac{\partial}{\partial \ov{z}_i})
-\sqrt{-1} \left( 2(\nabla ^{X} \nabla ^{X}J)
_{(\tfrac{\partial}{\partial  \ov{z}_i}, \tfrac{\partial}{\partial z_i})}  
+ [ R^{TX}(\tfrac{\partial}{\partial z_i}, 
\tfrac{\partial}{\partial \ov{z}_i}), J] \right)
\tfrac{\partial}{\partial  \ov{z}_l}, \tfrac{\partial}{\partial \ov{z}_m}
\right \rangle  d\ov{z}_l d\ov{z}_m\\
&\left.- \frac{1}{8\pi}\left(R^E +\tfrac{1}{2}\tr\left[R^{T^{(1,0)}X} \right]\right)
(\tfrac{\partial}{\partial  \ov{z}_l}, \tfrac{\partial}{\partial \ov{z}_m})
d\ov{z}_l d\ov{z}_m \right\} I_{\bC\otimes E}   .
\end{split}
\end{equation}

From  \eqref{g13}, \eqref{g35} and  \eqref{g11}, we have
\begin{equation}\label{g14}
\begin{split}
-&\left((L^0_2)^{-1}\underline{P}^{N^\bot}(\underline{\mQ}_2- \mO_2)
\underline{P}^N\right)(0,0)=\\
&= \left \{\frac{1}{4\pi} \tr\left[R^{T^{(1,0)}X}\right]
(\tfrac{\partial}{\partial z_i}, \tfrac{\partial}{\partial \ov{z}_i})
\right.
 - \frac{1}{2 \pi} \left\langle \Big(R^{TX}
(\tfrac{\partial}{\partial z_i}, \tfrac{\partial}{\partial \ov{z}_i})
- 2 \sqrt{-1} (\nabla ^{X} \nabla ^{X}J)_{(\tfrac{\partial}{\partial z_i},
 \tfrac{\partial}{\partial \ov{z}_i})} \Big)
\tfrac{\partial}{\partial z_l}, \tfrac{\partial}{\partial \ov{z}_l}
\right \rangle \\
&\left. \hspace{6mm}+\Big[\frac{1}{24\pi} \left\langle R^{TX}
(\tfrac{\partial}{\partial z_i}, \tfrac{\partial}{\partial \ov{z}_i})
\tfrac{\partial}{\partial  \ov{z}_l}, \tfrac{\partial}{\partial \ov{z}_m}
\right \rangle  
- \frac{1}{8\pi}\left(R^E +\tfrac{1}{2}\tr[R^{T^{(1,0)}X} ]\right)
(\tfrac{\partial}{\partial  \ov{z}_l}, \tfrac{\partial}{\partial \ov{z}_m})
\Big]d\ov{z}_l d\ov{z}_m \right\}  I_{\bC\otimes E}\\
&= \left \{\frac{3}{8\pi} \left\langle 
(\nabla ^X_{\tfrac{\partial}{\partial z_l}} J) \tfrac{\partial}{\partial z_j},
 (\nabla ^X_{\tfrac{\partial}{\partial \ov{z}_l}} J)
\tfrac{\partial}{\partial \ov{z}_j} \right\rangle 
- \frac{1}{12\pi}\left\langle R^{TX}
(\tfrac{\partial}{\partial \ov{z}_l}, \tfrac{\partial}{\partial \ov{z}_m})
\tfrac{\partial}{\partial z_i}, \tfrac{\partial}{\partial \ov{z}_i}
\right \rangle  d\ov{z}_l d\ov{z}_m \right.\\
&\left.\hspace{6mm}-\frac{1}{8\pi} 
R^E(\tfrac{\partial}{\partial  \ov{z}_l}, \tfrac{\partial}{\partial \ov{z}_m})
d\ov{z}_l d\ov{z}_m\right\}   I_{\bC\otimes E}.
\end{split}
\end{equation}
By \cite[(2.39)]{MM04a}, 
\begin{multline}\label{g60}
-(\cL_0^{-1}  P^{N^\bot} \mO_2  \underline{P}^N) (0,0)
=\frac{1}{2\pi} \left \{ \left\langle
R^{TX} (\tfrac{\partial}{\partial z_i},\tfrac{\partial}{\partial \ov{z}_j})
\tfrac{\partial}{\partial z_j}, 
\tfrac{\partial}{\partial \ov{z}_i}\right\rangle
+  R^E (\tfrac{\partial}{\partial z_i},
\tfrac{\partial}{\partial \ov{z}_i})\right\}  I_{\bC\otimes E}.
\end{multline}

By \eqref{g51}, \eqref{g52}, \eqref{0g9}, \eqref{0g9a}, 
\eqref{g8}, \eqref{g14},
 \eqref{g60}, and the discussion after \eqref{g52}, we get
\begin{equation}\label{g18}
\begin{split}
\pmb{b}_1(x) =& \frac{1}{\pi}\Big[\left\langle
R^{TX} (\tfrac{\partial}{\partial z_i},\tfrac{\partial}{\partial \ov{z}_j})
\tfrac{\partial}{\partial z_j}, 
\tfrac{\partial}{\partial \ov{z}_i}\right\rangle
+  R^E (\tfrac{\partial}{\partial z_j},
\tfrac{\partial}{\partial \ov{z}_j})
- \frac{1}{36}\sum_{kl} 
\Big|(\nabla ^X_{\tfrac{\partial}{\partial \ov{z}_k}} J)
\tfrac{\partial}{\partial \ov{z}_l} \Big|^2\Big]I_{\bC\otimes E}\\
&+ \frac{1}{36\pi} \left \langle
 (\nabla ^X_{\tfrac{\partial}{\partial \ov{z}_k}} J)
\tfrac{\partial}{\partial \ov{z}_l},
\tfrac{\partial}{\partial \ov{z}_m} \right \rangle 
\left \langle (\nabla ^X_{\tfrac{\partial}{\partial z_k}} J)
\tfrac{\partial}{\partial z_i},\tfrac{\partial}{\partial z_j} \right \rangle
 d\ov{z}_l \wedge d\ov{z}_m I_{\bC\otimes E} 
i_{\tfrac{\partial}{\partial \ov{z}_j}}\wedge 
i_{\tfrac{\partial}{\partial \ov{z}_i}}\\
&- \frac{1}{8\pi } \left(\frac{2}{3} 
\left \langle R^{TX} \tfrac{\partial}{\partial z_i}, 
\tfrac{\partial}{\partial \ov{z}_i}\right \rangle + R^E\right) 
(\tfrac{\partial}{\partial \ov{z}_l},\tfrac{\partial}{\partial \ov{z}_m})
 d\ov{z}_l\wedge d\ov{z}_m I_{\bC\otimes E} \\
&+ \frac{1}{2\pi } \left(\frac{2}{3} 
\left \langle R^{TX} \tfrac{\partial}{\partial z_i}, 
\tfrac{\partial}{\partial \ov{z}_i}\right \rangle
+ R^E\right) (\tfrac{\partial}{\partial z_l},\tfrac{\partial}{\partial z_m}) 
I_{\bC\otimes E} i_{\tfrac{\partial}{\partial \ov{z}_m}}
\wedge i_{\tfrac{\partial}{\partial \ov{z}_l}}.
\end{split}
\end{equation}
Moreover, we learn from \cite[Lemma 2.2]{MM04a} that 
\begin{equation}\label{g28}
r^X=-\left \langle  R^{TX} (e_i,e_j)e_i,e_j\right \rangle
= 8 \left \langle  R^{TX} (\tfrac{\partial}{\partial z_i},
\tfrac{\partial}{\partial \ov{z}_j})\tfrac{\partial}{\partial z_j},
\tfrac{\partial}{\partial \ov{z}_i} \right \rangle
-   \frac{1}{4}  |\nabla ^{X}J|^2.
\end{equation}
We are now ready to conclude. By \eqref{g18}  and \eqref{g28}, we get the formula \eqref{c1}
for $\pmb{b}_1$. 
We obtain then  \eqref{c2} by taking the trace of \eqref{c1}. 
The proof of Theorem \ref{t0.2} is complete.


\def\cprime{$'$} \def\cprime{$'$}
\providecommand{\bysame}{\leavevmode\hbox to3em{\hrulefill}\thinspace}
\providecommand{\MR}{\relax\ifhmode\unskip\space\fi MR }
\providecommand{\MRhref}[2]{%
  \href{http://www.ams.org/mathscinet-getitem?mr=#1}{#2}
}
\providecommand{\href}[2]{#2}

\end{document}